\documentclass[11pt]{article}
\usepackage{amsmath, amssymb}
\usepackage{graphicx}

\DeclareMathOperator{\dom}{dom}
\DeclareMathOperator{\gph}{gph}
\DeclareMathOperator{\clm}{clm}
\DeclareMathOperator{\Lipusc}{Lipusc}
\DeclareMathOperator{\Hof}{Hof}
\DeclareMathOperator{\uclm}{uclm}
\DeclareMathOperator{\extr}{extr}
\DeclareMathOperator{\spann}{span}

\DeclareMathOperator{\conv}{conv}

\DeclareMathOperator{\rank}{rank}

\newtheorem{theo}{Theorem}
\newtheorem{lem}{Lemma}
\newtheorem{prop}{Proposition}
\newtheorem{cor}{Corollary}
\newtheorem{rem}{Remark}
\newtheorem{exa}{Example}
\newtheorem{definiti}{Definition}

\newenvironment{dem}[1][Proof]{\noindent \textbf{#1.} }{\ \rule{0.5em}{0.5em}}
\begin{document}

\title{From calmness to Hoffman constants\ for linear semi-infinite
inequality systems\thanks{%
This research has been partially supported by Grant PGC2018-097960-B-C21
from MICINN, Spain, and ERDF, "A way to make Europe", European Union, and
Grant PROMETEO/2021/063 from Generalitat Valenciana, Spain.}}
\author{J. Camacho\thanks{%
Center of Operations Research, Miguel Hern\'{a}ndez University of Elche,
03202 Elche (Alicante), Spain (j.camacho@umh.es, canovas@umh.es,
parra@umh.es).} \and M.J. C\'{a}novas\footnotemark[2] \and J. Parra%
\footnotemark[2]}
\date{}
\maketitle

\begin{abstract}
In this paper we focus on different -global, semi-local and local- versions
of Hoffman type inequalities expressed in a variational form. In a first
stage our analysis is developed for generic multifunctions between metric
spaces and we finally deal with the feasible set mapping associated with
linear semi-infinite inequality systems (finitely many variables and
possibly infinitely many constraints) parameterized by their right-hand
side. The Hoffman modulus is shown to coincide with the supremum of
Lipschitz upper semicontinuity and calmness moduli when confined to
multifunctions with a convex graph and closed images in a reflexive Banach
space, which is the case of our feasible set mapping. Moreover, for this
particular multifunction a formula --only involving the system's left-hand
side-- of the global Hoffman constant is derived, providing a generalization
to our semi-infinite context of finite counterparts developed in the
literature. In the particular case of locally polyhedral systems, the paper
also provides a point-based formula for the (semi-local) Hoffman modulus in
terms of the calmness moduli at certain feasible points (extreme points when
the nominal feasible set contains no lines), yielding a practically
tractable expression for finite systems.

\textbf{Key words.} Hoffman constants, Lipschitz upper semicontinuity,
calmness, linear inequality systems, feasible set mapping.\newline

\bigskip

\noindent \textbf{Mathematics Subject Classification: } 90C31, 49J53, 90C34,
15A39, 90C05
\end{abstract}

\section{Introduction}

Concerning finite linear inequality systems parameterized by their
right-hand side, the celebrated Hoffman lemma \cite{Hof52} is a result of 
\emph{global} nature as far as it works for \emph{any} parameter making the
system consistent and \emph{any} point of the Euclidean space. We can also
find in the literature related \emph{semi-local} results as far as they work
around a nominal (given) parameter and any point in the Euclidean space,
leading to the concept of Hoffman constant at this parameter (see e.g. Az%
\'{e} and Corvellec \cite{AzCo02} and Z\u{a}linescu \cite{Za03}). In this
paper we relate these global and semi-local Hoffman constants with the \emph{%
local} concept of \emph{calmness modulus}, which involves parameters and
points, both around nominal ones. Our analysis is developed in a first step
in the context of generic multifunctions to move subsequently to the
particular case of the \emph{feasible set mapping} associated with a
parameterized linear semi-infinite inequality system%
\begin{equation}
\mathcal{\sigma }\left( b\right) :=\left\{ a_{t}^{\prime }x\leq
b_{t},\,\,\,t\in T\right\} ,  \label{eq_system}
\end{equation}%
where $T$ is a compact metric space, $t\mapsto a_{t}\in \mathbb{R}^{n}$ is a
fixed continuous function from $T$ to $\mathbb{R}^{n}$ and $b\equiv \left(
b_{t}\right) _{t\in T}\in C\left( T,\mathbb{R}\right) $ is the parameter to
be perturbed, $C\left( T,\mathbb{R}\right) $ being the space of continuous
functions from $T$ to $\mathbb{R}$. We are considering column-vectors and
the prime stands for transposition, so $x^{\prime }y$ denotes the usual
inner product of $x$ and $y$ in $\mathbb{R}^{n}.$ In this parametric
context, the feasible set mapping$,$ $\mathcal{F}:C\left( T,\mathbb{R}%
\right) \rightrightarrows \mathbb{R}^{n}$ is given by 
\begin{equation}
\mathcal{F}\left( b\right) :=\left\{ x\in \mathbb{R}^{n}\mid a_{t}^{\prime
}x\leq b_{t},\,\,\,t\in T\right\} .  \label{eq_feasible}
\end{equation}%
With respect to the topology, $\mathbb{R}^{n}$ is equipped with an arbitrary
norm, $\left\Vert \cdot \right\Vert ,$ with \emph{dual norm} given by $%
\left\Vert u\right\Vert _{\ast }=\max_{\left\Vert x\right\Vert \leq
1}\left\vert u^{\prime }x\right\vert ,$ and the parameter space $C\left( T,%
\mathbb{R}\right) $ is endowed with the supremum norm $\left\Vert
b\right\Vert _{\infty }:=\max_{t\in T}\left\vert b_{t}\right\vert .$

The particular case when $T$ is finite is included in this framework, in
which case $\mathcal{F}$ coincides with the polyhedral mapping considered in 
\cite{Hof52} and Hoffman lemma reads as the existence of some constant $%
\kappa \geq 0$ such that, for all $x\in \mathbb{R}^{n}$ and all $b\in \dom%
\mathcal{F}$ (the domain of $\mathcal{F)},$ 
\begin{equation}
d\left( x,\mathcal{F}\left( b\right) \right) \leq \kappa \max_{t\in T}\left[
a_{t}^{\prime }x-b_{t}\right] _{+}\text{,}  \label{eq_0002}
\end{equation}%
where $\left[ \alpha \right] _{+}:=\max \left\{ \alpha ,0\right\} $ is the
positive part of $\alpha \in \mathbb{R}.$ This result is of global nature as
far as it involves all points $x\in \mathbb{R}^{n}$ and all $b\in \dom%
\mathcal{F}.$ Since $\max_{t\in T}\left[ a_{t}^{\prime }x-b_{t}\right]
_{+}=d\left( b,\mathcal{F}^{-1}\left( x\right) \right) ,$ inequality (\ref%
{eq_0002}) can be written in a variational form as done in the following
paragraph for a generic multifunction.

Given a multifunction $\mathcal{M}:Y\rightrightarrows X$ between metric
spaces with both distances being denoted by $d,$ we say that the (global)
Hoffman property holds if there exists a constant $\kappa \geq 0$ such that 
\begin{equation}
d(x,\mathcal{M}(y))\leq \kappa d\left( y,\mathcal{M}^{-1}(x)\right) \text{
for all }x\in X\text{ and all }y\in \dom\mathcal{M},  \label{eq_hoffman}
\end{equation}%
where $d\left( x,\Omega \right) :=$ $\inf \left\{ d\left( x,\omega \right)
\mid \omega \in \Omega \right\} $ for $x\in X$ and $\Omega \subset X,$ with $%
\inf \emptyset :=+\infty ,$ so that $d\left( x,\emptyset \right) =+\infty .$
Since this paper is concerned with nonnegative constants, we use the
convention $\sup \emptyset :=0.$ Here $\dom\mathcal{M}$ is the domain of $%
\mathcal{M}$ (recall that $y\in \dom\mathcal{M\Leftrightarrow }$ $\mathcal{M}%
(y)\neq \emptyset )$ and $\mathcal{M}^{-1}$ denotes the inverse mapping of $%
\mathcal{M}$ (i.e. $y\in \mathcal{M}^{-1}\left( x\right) \Leftrightarrow
x\in \mathcal{M}(y)).$

Now we write a semi-local version of (\ref{eq_hoffman}) by fixing $y=%
\overline{y}.$ $\mathcal{M}$ is said to be \emph{Hoffman stable at }$%
\overline{y}\in \dom\mathcal{M}$ if there exists $\kappa \geq 0$ such that%
\begin{equation}
d(x,\mathcal{M}(\overline{y}))\leq \kappa d\left( \overline{y},\mathcal{M}%
^{-1}(x)\right) \text{ for all }x\in X\text{.}  \label{eq_Hof stable_b}
\end{equation}%
When the previous inequality (\ref{eq_Hof stable_b}) is only required to be
satisfied in a neighborhood of $\overline{x}\in \mathcal{M}(\overline{y})$
we are dealing with the \emph{calmness }of $\mathcal{M}$\ at $\left( 
\overline{y},\overline{x}\right) \in \gph\mathcal{M},$ the graph of $%
\mathcal{M}$. Formally, the calmness of $\mathcal{M}$\ at $\left( \overline{y%
},\overline{x}\right) \in \gph\mathcal{M},$ or equivalently the \emph{metric
subregularity} of $\mathcal{M}^{-1}$ \ at $\left( \overline{x},\overline{y}%
\right) $ (cf. \cite[Theorem 3H.3 and Exercise 3H.4]{DoRo}), is satisfied
when there exist a constant $\kappa \geq 0$ and a neighborhood $U$ of $%
\overline{x}$ such that 
\begin{equation}
d(x,\mathcal{M}(\overline{y}))\leq \kappa d\left( \overline{y},\mathcal{M}%
^{-1}(x)\right) \text{ for all }x\in U.  \label{eq subreg}
\end{equation}%
The infimum of constants $\kappa $ appearing in (\ref{eq_hoffman}), (\ref%
{eq_Hof stable_b}) and (\ref{eq subreg}) are called, respectively, the \emph{%
global Hoffman constant} of $\mathcal{M},$ the \emph{Hoffman modulus} of $%
\mathcal{M}$ at $\overline{y}\in \dom\mathcal{M},$ and the \emph{calmness
modulus }of $\mathcal{M}$ at $\left( \overline{y},\overline{x}\right) \in %
\gph\mathcal{M}.$ The three constants are denoted respectively by $\Hof%
\mathcal{M},$ $\Hof\mathcal{M}(\overline{y})$ and $\clm\mathcal{M}\left( 
\overline{y},\overline{x}\right) $ and, as a consequence of the definitions,
they may be written as follows: 
\begin{equation}
\begin{array}{l}
\Hof\mathcal{M}=\sup_{\left( y,x\right) \in \left( \dom\mathcal{M}\right) 
\mathcal{\times }X}\dfrac{d(x,\mathcal{M}(y))}{d\left( y,\mathcal{M}%
^{-1}(x)\right) },\medskip \\ 
\Hof\mathcal{M}(\overline{y})=\sup_{x\in X}\dfrac{d(x,\mathcal{M}(\overline{y%
}))}{d\left( \overline{y},\mathcal{M}^{-1}(x)\right) },\text{ }\overline{y}%
\in \dom\mathcal{M},\medskip \\ 
\clm\mathcal{M}\left( \overline{y},\overline{x}\right)
=\limsup_{x\rightarrow \overline{x}}\dfrac{d(x,\mathcal{M}(\overline{y}))}{%
d\left( \overline{y},\mathcal{M}^{-1}(x)\right) },\text{ }\left( \overline{y}%
,\overline{x}\right) \in \gph\mathcal{M},%
\end{array}
\label{eq_Hof_and_calm_modulus}
\end{equation}%
under the convention $\frac{0}{0}:=0,$ where $\lim \sup $ is understood as
the supremum (maximum, indeed) of all possible sequential upper limits
(i.e., with $(y,x)$ being replaced with elements of sequences $\left\{
(y_{r},x_{r})\right\} _{r\in \mathbb{N}}$ converging to $\left( \overline{y},%
\overline{x}\right) $ as $r\rightarrow \infty $).

Now we describe the main contributions of the paper. Clearly%
\begin{equation*}
\Hof\mathcal{M}=\sup_{\overline{y}\in \dom\mathcal{M}}\Hof\mathcal{M}(%
\overline{y}),
\end{equation*}%
and we wonder if a similar relationship between $\Hof\mathcal{M}(\overline{y}%
)$ and the supremum of all calmness moduli $\clm\mathcal{M}\left( \overline{y%
},x\right) $, with $x\in \mathcal{M}(\overline{y})$, works. Section 3 is
devoted to this question and Theorem \ref{The_calmness_Hoffman_convex} gives
a positive answer when $\gph\mathcal{M}$ is convex and $\mathcal{M}(%
\overline{y})$ is closed, $Y$ being a normed space and $X$ being a reflexive
Banach space. $\ $Some examples show that the convexity assumption is not
superfluous. Moreover, some intermediate constants as the \emph{Lipschitz
upper semicontinuity modulus }are also considered.

With respect to mapping $\mathcal{F}$ our focus is on formulae only
involving the system's coefficients for $\Hof\mathcal{F}$ and $\Hof\mathcal{F%
}(\overline{b}),$ which are established in Theorems \ref{Th_hof:global} and %
\ref{thm 1}, respectively. The first one extends to the current
semi-infinite framework previous results on finite linear systems (see,
e.g., Burke and Tseng \cite[Theorem 8]{BuTseng96}, Klatte and Thiere \cite[%
Theorem 2.7]{KlTh95}, Pe\~{n}a \emph{et al. }\cite[Formula (3)]{PVZ20}); for
comparative purposes, some details are gathered in Section 2. Theorem \ref%
{thm 1} provides a formula for $\Hof\mathcal{F}(\overline{b})$ in terms of
the $a_{t}$'s, the $\overline{b}_{t}$'s and some feasible points in the case
when our system $\sigma \left( \overline{b}\right) $ is for \emph{locally
polyhedral. }Specifically, from the referred Theorem \ref%
{The_calmness_Hoffman_convex}, we have that 
\begin{equation*}
\Hof\mathcal{F}(\overline{b})=\sup_{x\in \mathcal{F}(\overline{b})}\clm%
\mathcal{F}\left( \overline{b},x\right) ,
\end{equation*}%
and Theorem \ref{thm 1} refines this expression by reducing the supremum to
a smaller set (which turns out to be finite when $T$ also is). Then, making
use of the expression for $\clm\mathcal{F}\left( \overline{b},\overline{x}%
\right) $ established in Li \emph{et al. }\cite{LMY18} (recalled in Theorem %
\ref{Th_LiMengYan}), we derive the announced point-based formula for $\Hof%
\mathcal{F}(\overline{b}).$ Here we use the term `point-based' to emphasize
the fact that the expression for $\Hof\mathcal{F}(\overline{b})$ does not
involve parameters different from $\overline{b}$ or points outside $\mathcal{%
F}(\overline{b}).$ An alternative expression for $\Hof\mathcal{F}(\overline{b%
})$ appealing to points outside $\mathcal{F}(\overline{b})$ is given in \cite%
[Theorem 2.6]{AzCo02} (recalled in Theorem \ref{Th_Aze_Cor_subdif}). We
point out the fact that Theorem \ref{thm 1} yields a particularly tractable
procedure for computing $\Hof\mathcal{F}(\overline{b})$ when $T$ is finite.

In summary, the structure of the paper is as follows: Section 2 introduces
the necessary notation and gathers some preliminary results. Section 3
analyzes the relationships among different semi-local versions of Hoffman
and Lipschitz type properties for generic multifunctions and their moduli
(Lipschitz type properties are widely analyzed in the monographs \cite{DoRo,
KlKu02, mor06a, rw}). Section 3 also provides illustrative counter-examples.
Section 4 is focused on $\Hof\mathcal{F}$ and $\Hof\mathcal{F}(\overline{b})$%
, the latter in the case of locally polyhedral systems. Before establishing
the announced formula for $\Hof\mathcal{F}(\overline{b})$ some technical
geometrical results are proved. The paper finishes with a short section of
conclusions and perspectives.

\section{Preliminaries}

Given $S\subset \mathbb{R}^{k},$ $k\in \mathbb{N},$ we denote by \textrm{conv%
}$S,$ \textrm{cone}$S$ and $\mathrm{span}S$ the \emph{convex hull,} the 
\emph{conical convex hull and the linear hull }of $S$, respectively. It is
assumed that \textrm{cone}$S$ always contains the zero-vector $0_{k}$, in
particular \textrm{cone}$(\emptyset )=\{0_{k}\}.$ Moreover, $S^{\circ }$
denotes the (negative) \emph{polar }of $S$ given by 
\begin{equation*}
S^{\circ }:=\left\{ u\in \mathbb{R}^{k}\mid u^{\prime }x\leq 0,\text{ for
all }x\in S\right\}
\end{equation*}%
($S^{\circ }=\mathbb{R}^{k}$ if $S=\emptyset $). From the topological side, $%
\mathrm{int}S,$ $\mathrm{cl}S$ and $\mathrm{bd}S$ stand, respectively, for
the (topological) interior, closure, and boundary of $S.$ For a nonempty
convex set $C\subset \mathbb{R}^{k},$ $O^{+}C$ denotes its $\emph{recession}$
$\emph{cone}$ given by 
\begin{equation*}
O^{+}C:=\left\{ d\in \mathbb{R}^{k}\mid u+\alpha d\in C,\text{ for all }u\in
C\text{ and all }\alpha \geq 0\right\} ,
\end{equation*}%
while $\mathrm{end\,}C$ denotes its end set (introduced in \cite{HU05})\emph{%
\ }defined as 
\begin{equation*}
\mathrm{end\,}C:=\left\{ u\in \mathrm{cl\,}C\mid \nexists \mu >1\text{ such
that }\mu u\in \mathrm{cl\,}C\right\} .
\end{equation*}%
Moreover, $\mathrm{extr}C$ stands for the set of extreme points of $C$.
Recall that $x\in \mathrm{extr}C$ if $x\in C$ and it cannot be expressed as
a convex combination of two points of $C\backslash \{x\}.$ In any metric
space $\left( Z,d\right) ,$ the \emph{closed} ball centered at $z\in Z$ with
radius $r>0$ is denoted by $B\left( z,r\right) ,$ whereas $B\left(
S,r\right) :=\left\{ z\in Z\mid d\left( z,S\right) \leq r\right\} ,$ for $%
S\subset Z,$ denotes the $r$\emph{-enlargement} of $S.$

For comparative purposes, the next theorem gathers some results in the
literature on $\Hof\mathcal{F}$ when confined to finite linear systems,
where $C\left( T,\mathbb{R}\right) \equiv \mathbb{R}^{m}$ for some $m\in 
\mathbb{N}.$ It is adapted to our current notation and to our choice of
norms. The first two expressions come from \cite[Formulae (3) and (4)]{PVZ20}
(see also \cite[Theorem 2.7]{KlTh95} when $\mathbb{R}^{n}$ is endowed with
the Euclidean norm), while the third one can be derived from \cite[Theorem 8]%
{BuTseng96}, where a dual approach is followed. The last one appeals to the
set 
\begin{equation*}
W_{2}:=\left\{ y\in \mathbb{R}_{+}^{m}\mid \left\{ a_{t},\text{ }t\in 
\mathrm{supp}\left( y\right) \right\} \text{ lin. indep.}\right\} ,
\end{equation*}%
where $\mathbb{R}_{+}^{m}$ is formed by the vectors of $\mathbb{R}^{m}$
having non-negative coordinates and $\mathrm{supp}\left( y\right) :=\left\{
t\in \{1,...,m\}\mid y_{t}\neq 0\right\} $ is the support of $y;$ indeed $%
W_{2}$ is considered as a subset of the dual space of $\mathbb{R}^{m},$
which we are identifying with $\mathbb{R}^{m}$ itself$.$

\begin{theo}
\label{TH_PVZ_const}Consider the feasible set mapping $\mathcal{F}$ defined
in (\ref{eq_feasible}) and assume that $T$ is finite. We have%
\begin{eqnarray}
\Hof\mathcal{F} &=&\max_{_{\substack{ J\subset T  \\ 0_{n}\notin \conv%
\left\{ a_{t},~t\in J\right\} }}}d_{\ast }\left( 0_{n},\conv\left\{
a_{t},~t\in J\right\} \right) ^{-1}  \label{eq_1} \\
&=&\max_{\substack{ J\subset T,\text{ }\rank A_{J}=\rank A  \\ \left\{ a_{t},%
\text{ }t\in J\right\} \text{ lin. indep.}}}d_{\ast }\left( 0_{n},\conv%
\left\{ a_{t},~t\in J\right\} \right) ^{-1}  \label{eq_2} \\
&=&\sup \left\{ \left\Vert y\right\Vert _{1}\mid y\in W_{2},\text{ }%
\left\Vert A^{\prime }y\right\Vert _{\ast }=1\right\} ,  \label{eq_3}
\end{eqnarray}%
where $A_{J}$ and $A$ stand for the matrices whose rows are $a_{t}^{\prime
}, $ with $t\in J$ and $t\in T,$ respectively, and $d_{\ast }$ stands for
the distance associated with the dual norm $\left\Vert \cdot \right\Vert
_{\ast }.$
\end{theo}

\begin{dem}
According to \cite[Formula (3)]{PVZ20} and the subsequent comments therein,
to establish (\ref{eq_1}) we only have to prove that condition $0_{n}\notin %
\conv\left\{ a_{t},~t\in J\right\} $ is equivalent to the consistency of
system $\left\{ a_{t}^{\prime }x<0,~t\in J\right\} ,$ and this follows, for
instance, from equivalence $\left( iv\right) \Leftrightarrow \left( v\right) 
$ in \cite[Theorem 6.1]{libro}. Equality (\ref{eq_2}) comes from \cite[%
Formula (4)]{PVZ20} with the trivial observation that instead of all
linearly independent $\left\{ a_{t},\text{ }t\in J\right\} ,$ with $J\subset
T$, we can confine ourselves to those which are maximal with respect to the
inclusion order. Indeed, the result also follows from (\ref{eq_1}), since
the sufficiency of considering those $\left\{ a_{t},\text{ }t\in J\right\} $
which are linearly independent comes from \cite[Lemma 3.1]{AzCo02}.

Formula (\ref{eq_3}) comes from \cite[Theorem 8]{BuTseng96}. Let us comment
that we can, alternatively, see the relationship between the second and the
third expression by observing that, for any $y\in \mathbb{R}_{+}^{m},$ $%
y\neq 0_{m},$ 
\begin{equation*}
\frac{1}{\left\Vert y\right\Vert _{1}}A^{\prime }y=\frac{1}{\left\Vert
y\right\Vert_{1}}\sum\nolimits_{i=1}^{m}y_{t}a_{t}\in \conv\left\{
a_{t},~t\in \mathrm{supp}\left( y\right) \right\} ,
\end{equation*}%
and that $\left\Vert A^{\prime }y\right\Vert _{\ast }=1$ is equivalent to $%
\left\Vert y\right\Vert _{1}=\left\Vert \frac{1}{\left\Vert y\right\Vert _{1}%
}A^{\prime }y\right\Vert _{\ast }^{-1}.$
\end{dem}

Generalizations of Hoffman constants to infinite dimensional spaces or to
convex functions playing the role of the distance function can be found in 
\cite{BuTseng96}. Many other authors have contributed to the study of
Hoffman constants and their relationship with other concepts (as Lipschitz
constants). Additional references can be obtained from the reference list of
the papers above mentioned as well as \cite{AzCo02} and \cite{Za03}, among
others. At this moment we also cite Belousov and Andronov \cite{BeAn99}, Li 
\cite{Li94} and Robinson \cite{Rob81}.

The following theorem provides formulae for $\Hof\mathcal{F}(\overline{b}),$
with $\overline{b}\in \dom\mathcal{F},$ and $\clm\mathcal{F}\left( \overline{%
b},\overline{x}\right) ,$ with $\left( \overline{b},\overline{x}\right) \in %
\gph\mathcal{F}$ through points outside $\mathcal{F}(\overline{b}).$ They
appeal to the supremum function $f_{b}:\mathbb{R}^{n}\rightarrow \mathbb{R},$
with $b\in C\left( T,\mathbb{R}\right) ,$ given by 
\begin{equation*}
f_{b}\left( x\right) :=\sup_{t\in T}\left( a_{t}^{\prime }x-b_{t}\right) ,%
\text{ for }x\in \mathbb{R}^{n},
\end{equation*}%
which is known to be convex on $\mathbb{R}^{n}$. For each $x\in \mathbb{R}%
^{n},$ we consider the subset of indices 
\begin{equation*}
J_{b}\left( x\right) =\left\{ t\in T\mid a_{t}^{\prime }x-b_{t}=f_{b}\left(
x\right) \right\} .
\end{equation*}%
The well-known Valadier's formula works by virtue of the Ioffe-Tikhomirov
theorem (see e.g. \cite[Theorem 2.4.18]{ZA02}), yielding 
\begin{equation*}
\partial f_{b}\left( x\right) =\mathrm{conv}\left\{ a_{t},\text{ }t\in
J_{b}\left( x\right) \right\} ,
\end{equation*}%
where $\partial f_{b}\left( x\right) $ stands for the usual subdifferential
of convex analysis (see\ e.g. \cite{Rock70}).

\begin{theo}
\label{Th_Aze_Cor_subdif} The following statements hold:

$\left( i\right) $ \emph{\cite[Theorem 2.6]{AzCo02}} For any $\overline{b}%
\in \dom\mathcal{F},$ one has 
\begin{eqnarray*}
\Hof\mathcal{F}(\overline{b}) &=&\sup_{f_{\overline{b}}\left( x\right)
>0}d_{\ast }\left( 0_{n},\partial f_{\overline{b}}\left( x\right) \right)
^{-1} \\
&=&\sup_{f_{\overline{b}}\left( x\right) >0}d_{\ast }\left( 0_{n},\mathrm{%
conv}\left\{ a_{t},\text{ }t\in J_{\overline{b}}\left( x\right) \right\}
\right) ^{-1};
\end{eqnarray*}

$\left( ii\right) $ \emph{\cite[Theorem 1]{KNT10} }For any $\left( \overline{%
b},\overline{x}\right) \in \gph\mathcal{F},$ 
\begin{eqnarray*}
\mathrm{clm}\mathcal{F}(\overline{b},\overline{x}) &=&\underset{x\rightarrow 
\overline{x},~f_{\overline{b}}\left( x\right) >0}{\lim \sup }d_{\ast }\left(
0_{n},\partial f_{\overline{b}}\left( x\right) \right) ^{-1} \\
&=&\underset{x\rightarrow \overline{x},~f_{\overline{b}}\left( x\right) >0}{%
\lim \sup }d_{\ast }\left( 0_{n},\mathrm{conv}\left\{ a_{t},\text{ }t\in J_{%
\overline{b}}\left( x\right) \right\} \right) ^{-1}.
\end{eqnarray*}
\end{theo}

\begin{rem}
\label{Rem1}\emph{Observe that }$\overline{b}\in \dom\mathcal{F}$ \emph{and }%
$f_{\overline{b}}\left( x\right) >0$\emph{\ mean that }$\sigma \left( 
\overline{b}\right) $ \emph{is consistent (it has some feasible solution)
but }$x\notin \mathcal{F}(\overline{b})$\emph{); in this case}$,$ $%
0_{n}\notin \mathrm{conv}\{a_{t},~t\in J_{\overline{b}}\left( x\right) \},$ 
\emph{since }$x$ \emph{is not a global minimizer of the convex function }$f_{%
\overline{b}}.$ \emph{Actually,} \emph{\cite[Theorem 2.6]{AzCo02}} \emph{\
is formulated in terms of} $(\Hof\mathcal{F}(\overline{b}))^{-1}$\emph{,
which is called there the} condition number of $f_{\overline{b}}$ at level $%
0 $; \emph{in the terminology of \cite{KNT10}, observe that }$(\mathrm{clm}%
\mathcal{F}(\overline{b},\overline{x}))^{-1}$ \emph{is the} error bound
modulus \emph{(also known as} conditioning rate \emph{\cite{Penot10})} of $%
f_{\overline{b}}$ at $\overline{x}.$
\end{rem}

The following theorem is devoted to the computation of $\mathrm{clm}\mathcal{%
F}(\overline{b},\overline{x}),$ $(\overline{b},\overline{x})\in \gph\mathcal{%
F},$ through a point-based formula (expressed exclusively in terms of the
system's coefficients and the nominal point $\overline{x}).$ Now we
introduce some extra notation. Given a fixed $\overline{b}\in \dom\mathcal{F}
$, for any $x\in $ $\mathcal{F}(\overline{b}),$ we consider (for simplicity,
since there will be no ambiguity, we omit the dependence on $\overline{b}$) 
\begin{equation*}
T\left( x\right) :=\left\{ t\in T\mid a_{t}^{\prime }x-\overline{b}%
_{t}=0\right\} ,
\end{equation*}%
the \emph{subset of active indices }of system $\sigma \left( \overline{b}%
\right) $ at\emph{\ }$x;$ i.e., $T\left( x\right) =J_{\overline{b}}\left(
x\right) $ if $f_{\overline{b}}\left( x\right) =0,$ while $T\left( x\right)
=\emptyset $ if $f_{\overline{b}}\left( x\right) <0$ (i.e., if $x$ is a
strict solution --\emph{Slater point}-- of the system). Let $A\left(
x\right) $ be the corresponding \emph{active cone} at $x;$ i.e., 
\begin{equation*}
A\left( x\right) :=\mathrm{cone}\left\{ a_{t},\text{ }t\in T\left( x\right)
\right\} 
\end{equation*}%
(recall that $A\left( x\right) =\{0_{n}\}$ if $T\left( x\right) =\emptyset ).
$ We also consider the family $\mathcal{D}\left( x\right) $ of subsets $%
D\subset T\left( x\right) $ such that system 
\begin{equation}
\left\{ 
\begin{tabular}{rl}
$a_{t}^{\prime }d=1,$ & $t\in D,$ \\ 
$a_{t}^{\prime }d<1,$ & $t\in T\left( x\right) \setminus D$%
\end{tabular}%
\right\}   \label{SystemD2}
\end{equation}%
is consistent (in the variable $d\in \mathbb{R}^{n});$ i.e., $\{a_{t},$ $%
t\in D\}$ is contained in some hyperplane which leaves $\{0_{n}\}\cup
\{a_{t},t\in T\left( x\right) \setminus D\}$ on one of its two associated
open half-spaces. With this notation, the next theorem generalizes the
corresponding finite version established in \cite[Theorem 4]{CLPT14}. It
appeals to the following \emph{regularity condition} at $\overline{x}$:
\textquotedblleft There exists a neighborhood $W$ of $\overline{x}$ such
that 
\begin{equation}
\mathcal{F}(\overline{b})\cap W=\left( \overline{x}+A\left( \overline{x}%
\right) ^{\circ }\right) \cap W.\text{\textquotedblright }
\label{eq_regularity_condition}
\end{equation}%
Observe that this condition is held at all points of polyhedral sets and,
for instance, at the vertex of the ice-cream cone.

\begin{theo}
\emph{\cite[Corollary 2.1, Remark 2.3 and Corollary 3.2]{LMY18}} \label%
{TH_calmness_mod}\label{Th_LiMengYan} Let $\overline{x}\in \mathcal{F}(%
\overline{b})$ such that $f_{\overline{b}}\left( \overline{x}\right) =0$ and
assume that the regularity condition (\ref{eq_regularity_condition}) is held
at $\overline{x}$. Then%
\begin{equation}
\clm\mathcal{F}\left( \overline{b},\overline{x}\right) =d_{\ast }\left(
0_{n},\mathrm{end}\partial f_{\overline{b}}\left( \overline{x}\right)
\right) ^{-1}=\sup_{D\in \mathcal{D}\left( \overline{x}\right) }d_{\ast
}\left( 0_{n},\mathrm{conv}\left\{ a_{t},\text{ }t\in D\right\} \right)
^{-1}.  \label{eq_calmness_modulus_LMY}
\end{equation}
\end{theo}

\begin{rem}
\label{Rem_end}\emph{Although condition (\ref{eq_regularity_condition}) is
not superfluous for establishing the first equality in (\ref%
{eq_calmness_modulus_LMY}) as \cite[Example 3.3]{LMY18} shows (see also
Example \ref{Exa1}), the second equality does work for semi-infinite systems
(\ref{eq_system}) without any additional condition. Indeed, from \cite[%
Corollary 2.1 and Remark 2.3]{LMY18} we can deduce} 
\begin{equation}
\cup _{D\in \mathcal{D}\left( \overline{x}\right) }\mathrm{conv}\left\{
a_{t},\text{ }t\in D\right\} \subset \mathrm{end}\partial f_{\overline{b}%
}\left( \overline{x}\right) \subset \mathrm{cl}\left( \cup _{D\in \mathcal{D}%
\left( \overline{x}\right) }\mathrm{conv}\left\{ a_{t},\text{ }t\in
D\right\} \right) .  \label{eq_end}
\end{equation}
\end{rem}

\section{From calmness to Hoffman constants for a ge\-neric multifunction}

The purpose of this section is to analyze the relationship among different
Hoffman and Lipschitz type properties, including the known Lipschitz upper
semicontinuity that goes back to the classical work of Robinson \cite{Rob81}%
. At the beginning of this section $\mathcal{M}:Y\rightrightarrows X$ is a
generic multifunction between metric spaces $Y$ and $X$. Later we will need
further structure. To start with, observe that alternatively to (\ref{eq_Hof
stable_b}) we can write the Hoffman stability of $\mathcal{M}$\emph{\ }at%
\emph{\ }$\overline{y}\in \dom\mathcal{M}$ in terms of the existence of $%
\kappa \geq 0$ such that%
\begin{equation*}
d(x,\mathcal{M}(\overline{y}))\leq \kappa d\left( y,\overline{y}\right) 
\text{ for all }\left( y,x\right) \in \gph\mathcal{M},
\end{equation*}%
while the calmness of $\mathcal{M}$ at $\left( \overline{y},\overline{x}%
\right) \in \gph\mathcal{M}$, introduced in (\ref{eq subreg}) in terms of
the (equivalent) metric subregularity of $\mathcal{M}^{-1},$ writes as the
existence of neighborhoods $V$ of $\overline{y}$ and $U$ of $\overline{x}$
along with a constant $\kappa \geq 0$ such that 
\begin{equation*}
d(x,\mathcal{M}(\overline{y}))\leq \kappa d\left( y,\overline{y}\right) 
\text{ for all }\left( y,x\right) \in \left( V\times U\right) \mathcal{\cap }%
\gph\mathcal{M}.
\end{equation*}%
Moreover, the following equalities constitute well-known alternative
expressions to (\ref{eq_Hof_and_calm_modulus}) for the corresponding moduli 
\begin{eqnarray}
\Hof\mathcal{M}(\overline{y}) &=&\sup_{\left( y,x\right) \in \gph\mathcal{M}}%
\frac{d(x,\mathcal{M}(\overline{y}))}{d\left( y,\overline{y}\right) },
\label{eq_hof_calm_def} \\
\clm\mathcal{M}\left( \overline{y},\overline{x}\right) &=&\limsup 
_{\substack{ (y,x)\rightarrow \left( \overline{y},\overline{x}\right)  \\ %
\left( y,x\right) \in \gph\mathcal{M}}}\frac{d(x,\mathcal{M}(\overline{y}))}{%
d\left( y,\overline{y}\right) }.  \notag
\end{eqnarray}

Recall that $\mathcal{M}$ is said to be \emph{Lipschitz upper semicontinuous 
}at $\overline{y}\in \dom\mathcal{M}$ if there exists a neighborhood $V$ of $%
\overline{y}$ along with a constant $\kappa \geq 0$ such that 
\begin{equation}
d(x,\mathcal{M}(\overline{y}))\leq \kappa d\left( y,\overline{y}\right) 
\text{ for all }\left( y,x\right) \in \left( V\times \mathbb{R}^{n}\right) 
\mathcal{\cap }\gph\mathcal{M}.  \label{eq Lipschitz_usc}
\end{equation}%
Here we borrow the terminology from \cite{KlKu02} or \cite{Uderzo21},
although this property, introduced in \cite{Rob81} as\textbf{\ }\emph{upper
Lipschitz continuity}\textbf{, }has been also popularized as\emph{\ outer
Lipschitz continuity }(see \cite{DoRo}). Equivalently, (\ref{eq
Lipschitz_usc}) may be written as $e(\mathcal{M}(y),\mathcal{M}(\overline{y}%
))\leq \kappa d\left( y,\overline{y}\right) $ for all $y\in V,$ where $%
e\left( A,B\right) :=\sup_{x\in A}d\left( x,B\right) $ is the \emph{%
Hausdorff excess} of $A$ over $B,$ with $A,B\subset X.$ The associated \emph{%
Lipschitz upper semicontinuity modulus}, denoted by $\Lipusc\mathcal{M}(%
\overline{y}),$ is defined as the infimum of constants $\kappa $ satisfying (%
\ref{eq Lipschitz_usc}) for some associated $V.$

In the next definition, given $\overline{y}\in \dom\mathcal{M}$ and $%
\varepsilon >0$, the mapping $\mathcal{M}_{\varepsilon }:Y\rightrightarrows
X $ is defined by%
\begin{equation*}
\mathcal{M}_{\varepsilon }\left( y\right) :=\mathcal{M}\left( y\right) \cap
B\left( \mathcal{M}\left( \overline{y}\right) ,\varepsilon \right) \text{
for }y\in Y.
\end{equation*}%
(For simplicity in the notation we obviate the dependence of $\mathcal{M}%
_{\varepsilon }$ on $\overline{y}.$)

\begin{definiti}
\label{Def uclm Hof}Given $\overline{y}\in \dom\mathcal{M},$ we say that $%
\mathcal{M}$ is \emph{uniformly calm} at $\overline{y}$ if there exist a
neighborhood $V$ of $\overline{y}$ along with $\varepsilon >0$ and $\kappa
\geq 0$ such that 
\begin{equation}
d(x,\mathcal{M}(\overline{y}))\leq \kappa d\left( y,\overline{y}\right) 
\text{ for all }y\in V\text{ and all }x\in \mathcal{M}_{\varepsilon }\left(
y\right) ,  \label{eq uniform calmness}
\end{equation}%
or, equivalently, if $\mathcal{M}_{\varepsilon }$ is Lipschitz upper
semicontinuous at $\overline{y}$ for some $\varepsilon >0.$
\end{definiti}

The corresponding modulus naturally appear. Specifically, we call \emph{%
modulus of uniform calmness of }$\mathcal{M}$ at $\overline{y},$ denoted by $%
\uclm\mathcal{M}(\overline{y}),$ to the infimum of constants $\kappa $
satisfying (\ref{eq uniform calmness}) for some associated $V$ and $%
\varepsilon >0$. It is straightforward to check that 
\begin{equation}
\uclm\mathcal{M}(\overline{y})=\inf_{\varepsilon >0}\Lipusc\mathcal{M}%
_{\varepsilon }(\overline{y}).  \label{eq uclm as inf eps}
\end{equation}%
Roughly speaking, the uniform calmness of $\mathcal{M}$ at $\overline{y}$
entails the calmness of $\mathcal{M}$ at any $(\overline{y},x)$ for all $%
x\in \mathcal{M}\left( \overline{y}\right) $ with the same calmness constant 
$\kappa ,$ the same neighborhood $V$ of $\overline{y},$ and a common radius $%
\varepsilon $ for all neighborhoods of points $x\in \mathcal{M}\left( 
\overline{y}\right) ,$ say $U_{x}:=B(x,\varepsilon )$. Example \ref{Exa_1}
below shows that the calmness of $\mathcal{M}$ at $(\overline{y},x)$ for all 
$x\in \mathcal{M}\left( \overline{y}\right) $ does not ensure the uniform
calmness of $\mathcal{M}$ at $\overline{y}.$

As it occurs with the calmness property, the uniform calmness turns out to
be equivalent to a certain metric regularity type property, showing that
neighborhood $V$ in Definition \ref{Def uclm Hof} is redundant. The key fact
is that points $x\in \mathcal{M}(y)$ which are required to satisfy (\ref{eq
uniform calmness}) are those which are sufficiently close to $\mathcal{M}%
\left( \overline{y}\right) .$ This comment, which was already pointed out
for polyhedral multifunctions in \cite{Rob81} (see the corollary after
Proposition 1 therein), is formalized in the following proposition.

\begin{prop}
\label{Prop eqiv k>0}Let $\overline{y}\in \dom\mathcal{M}.$ For any $\kappa
>0,$ the following conditions are equivalent:

$\left( i\right) $ There exist a neighborhood $V$ of $\overline{y}$ and $%
\varepsilon >0$ such that (\ref{eq uniform calmness}) holds;

$\left( ii\right) $ There exists $\varepsilon >0$ such that (\ref{eq_Hof
stable_b}) holds when restricted to those $x\in B\left( \mathcal{M}\left( 
\overline{y}\right) ,\varepsilon \right) .$
\end{prop}

\begin{dem}
Let us establish the nontrivial implication `$\left( i\right) \Rightarrow
\left( ii\right) $'. Consider $V$ and $\varepsilon $ as in statement $\left(
i\right) .$ Take $\varepsilon _{1}>0$ such that $B(\overline{y},\varepsilon
_{1})\mathbb{\subset }V$ and define $\varepsilon _{2}:=\min \{\varepsilon
,\kappa \varepsilon _{1}\}>0.$ Let us see that $\left( ii\right) $ holds for 
$\varepsilon _{2}>0.$ Take $x\in B\left( \mathcal{M}\left( \overline{y}%
\right) ,\varepsilon \right) $ and consider $y\in \mathcal{M}^{-1}\left(
x\right) .$ Now, we distinguish between two cases:

If $d\left( y,\overline{y}\right) \leq \varepsilon _{1},$ then $y\in V$ and,
since we also have $x\in B\left( \mathcal{M}\left( \overline{y}\right)
,\varepsilon \right) $ (recall that $\varepsilon _{2}\leq \varepsilon $),
from $\left( i\right) $ we conclude the aimed inequality $d(x,\mathcal{M}(%
\overline{y}))\leq \kappa d\left( y,\overline{y}\right) .$

Otherwise, if $d\left( y,\overline{y}\right) \geq \varepsilon _{1},$ then $%
d(x,\mathcal{M}(\overline{y}))\leq \varepsilon _{2}\leq \kappa \varepsilon
_{1}\leq \kappa d\left( y,\overline{y}\right) .$
\end{dem}

\begin{rem}
\emph{The statement of Proposition \ref{Prop eqiv k>0} does not hold for }$%
\kappa =0.$ \emph{To see this, take }$\mathcal{M}:\mathbb{R}\longrightarrow 
\mathbb{R}$ \emph{(single-valued) given by }$\mathcal{M}\left( y\right)
:=\max \left\{ 0,y-1\right\} $\emph{\ and let }$\overline{y}=0.$\emph{\
Clearly }$\left( i\right) $\emph{\ holds for} $V=\left] -1,1\right[ $ \emph{%
and} $\kappa =0,$\emph{\ whereas }$\left( ii\right) $ \emph{works for} $%
\varepsilon >0$\emph{\ if and only if }$\kappa \geq \varepsilon /\left(
1+\varepsilon \right) .$
\end{rem}

\begin{cor}
\label{Cor_unif_calm}Let $\overline{y}\in \dom\mathcal{M}.$ We have:

$\left( i\right) $ $\mathcal{M}$ is uniformly calm at $\overline{y}$ if and
only if there exist $\varepsilon >0$ and $\kappa \geq 0$ such that 
\begin{equation}
d(x,\mathcal{M}(\overline{y}))\leq \kappa d\left( \overline{y},\mathcal{M}%
^{-1}(x)\right) \text{ for all }x\in B\left( \mathcal{M}\left( \overline{y}%
\right) ,\varepsilon \right) .  \label{eq_002}
\end{equation}

$\left( ii\right) $ The modulus of uniform calmness can be expressed as
follows%
\begin{equation*}
\uclm\mathcal{M}(\overline{y})=\inf \left\{ \kappa \geq 0\mid \exists
\varepsilon >0\text{ such that }(\ref{eq_002})\text{ holds}\right\} .
\end{equation*}
\end{cor}

\begin{dem}
Both $\left( i\right) $ and $\left( ii\right) $ come from the fact that
uniform calmness at $\overline{y}$ with associated elements $V,$ $%
\varepsilon >0$ and $\kappa \geq 0$ in (\ref{eq uniform calmness}) entails
the same property with $V,$ $\varepsilon >0$ and $\widetilde{\kappa }>\kappa
.$ Hence the conclusions follow straightforwardly from Proposition \ref{Prop
eqiv k>0}.
\end{dem}

Next we provide characterizations of $\Lipusc\mathcal{M}(\overline{y})$ and $%
\uclm\mathcal{M}(\overline{y})$ in terms of certain upper limits, which
allow for a better understanding of these concepts and a clear relationship
among all moduli introduced in the paper.

\begin{prop}
Let $\mathcal{M}:Y\rightrightarrows X$ be a multifunction between metric
spaces and let $\overline{y}\in \dom\mathcal{M},\,$\ then

$\left( i\right) $ $\Lipusc\mathcal{M}(\overline{y})=\limsup\limits_{y%
\rightarrow \overline{y}}\left( \sup\limits_{x\in \mathcal{M}(y)}\dfrac{d(x,%
\mathcal{M}(\overline{y}))}{d\left( y,\overline{y}\right) }\right) ;$

$\left( ii\right) ~\uclm\mathcal{M}(\overline{y})=\limsup\limits_{d(x,%
\mathcal{M}(\overline{y}))\rightarrow 0}\dfrac{d(x,\mathcal{M}(\overline{y}))%
}{d\left( \overline{y},\mathcal{M}^{-1}(x)\right) }$.
\end{prop}

\begin{dem}
$\left( i\right) $ For the sake of simplicity, let us denote by $s$ the
right-hand side of $\left( i\right) $ and 
\begin{equation}
K:=\left\{ \kappa \geq 0\mid \exists V\text{ neighborhood of }\overline{y}%
\text{ verifying (\ref{eq Lipschitz_usc})}\right\} .  \label{eq K}
\end{equation}%
We start by establishing inequality `$\leq $'. Since $\Lipusc\mathcal{M}(%
\overline{y})=\inf K,$ we can write $\Lipusc\mathcal{M}(\overline{y}%
)=\lim_{r\rightarrow \infty }\kappa _{r}$ for some $\left\{ \kappa
_{r}\right\} \subset K.$ For each $r$ take a neighborhood $V_{r}$ associated
with $\kappa _{r}$ according to (\ref{eq K}) and define 
\begin{equation*}
\overline{\kappa }_{r}:=\sup_{y\in V_{r}\cap B\left( \overline{y},1/r\right)
}\left( \sup_{x\in \mathcal{M}(y)}\frac{d(x,\mathcal{M}(\overline{y}))}{%
d\left( y,\overline{y})\right) }\right) \leq \kappa _{r}.
\end{equation*}%
By definition $\overline{\kappa }_{r}\in K,$ having $V_{r}\cap B\left( 
\overline{y},1/r\right) $ as an associated neighborhood, so that we have $%
\Lipusc\mathcal{M}(\overline{y})=\lim_{r\rightarrow \infty }\overline{\kappa 
}_{r}.$

Finally, for each $r,$ consider any $y_{r}\in V_{r}\cap B\left( \overline{y}%
,1/r\right) $ such that $\overline{\kappa }_{r}-\frac{1}{r}\leq \sup_{x\in 
\mathcal{M}(y_{r})}\frac{d(x,\mathcal{M}(\overline{y}))}{d\left( y_{r},%
\overline{y})\right) }\leq \overline{\kappa }_{r}.$ Obviously, $\left\{
y_{r}\right\} _{r\in \mathbb{N}}$ converges to $\overline{y},$ and then 
\begin{equation*}
\Lipusc\mathcal{M}(\overline{y})=\lim_{r\rightarrow \infty }\sup_{x\in 
\mathcal{M}(y_{r})}\frac{d(x,\mathcal{M}(\overline{y}))}{d\left( y_{r},%
\overline{y})\right) }\leq s.
\end{equation*}

In order to prove `$\geq $' in $\left( i\right) ,$ we may assume the
nontrivial case $s>0$ and write 
\begin{equation*}
s=\lim_{r\rightarrow \infty }\sup_{x\in \mathcal{M}(\widetilde{y}_{r})}\frac{%
d(x,\mathcal{M}(\overline{y}))}{d\left( \widetilde{y}_{r},\overline{y}%
\right) },
\end{equation*}%
for some $\left\{ \widetilde{y}_{r}\right\} _{r\in \mathbb{N}}$ converging
to $\overline{y}.$ It is clear that we may replace $\left\{ \widetilde{y}%
_{r}\right\} _{r\in \mathbb{N}}$ with a suitable subsequence (denoted as the
whole sequence for simplicity) such that $\widetilde{y}_{r}\in V_{r}$, and
then 
\begin{equation*}
s\leq \lim_{r\rightarrow \infty }\kappa _{r}=\Lipusc\mathcal{M}(\overline{y}%
).
\end{equation*}

$\left( ii\right) $ The procedure is analogous to the previous one by
considering 
\begin{equation*}
\widehat{K}=\left\{ \kappa \geq 0\mid \exists \varepsilon >0\text{ such that 
}(\ref{eq_002})\text{ holds}\right\} .
\end{equation*}
\end{dem}

As a direct consequence of the expressions in (\ref{eq_hof_calm_def}) for $%
\clm\mathcal{M}\left( \overline{y},\overline{x}\right) $ and $\Hof\mathcal{M}%
\left( \overline{y}\right) ,$ together with (\ref{eq uclm as inf eps}) and
the previous proposition, we conclude the following corollary. Observe that
the smaller $\varepsilon >0,$ the smaller $\Lipusc\mathcal{M}_{\varepsilon
}\left( \overline{y}\right) ,$ and $\Lipusc\mathcal{M}\left( \overline{y}%
\right) $ corresponds to $\varepsilon =+\infty .$

\begin{cor}
\label{Cor three inequalities}Let $\overline{y}\in \dom\mathcal{M}.$ We have 
\begin{equation}
\sup_{x\in \mathcal{M}(\overline{y})}\clm\mathcal{M}\left( \overline{y}%
,x\right) \leq \uclm\mathcal{M}\left( \overline{y}\right) \leq \Lipusc%
\mathcal{M}\left( \overline{y}\right) \leq \Hof\mathcal{M}\left( \overline{y}%
\right) .  \label{eq_inequalities chain}
\end{equation}
\end{cor}

\begin{rem}
\label{Rem implication chain}\emph{The previous corollary yields }$\left(
i\right) \Rightarrow \left( ii\right) \Rightarrow \left( iii\right)
\Rightarrow \left( iv\right) ,$ \emph{where:}

$\left( i\right) $ $\mathcal{M}$ \emph{is Hoffman stable at} $\overline{y};$

$\left( ii\right) $ $\mathcal{M}$ \emph{is Lipschitz upper semicontinuous}$\ 
$\emph{at} $\overline{y};$

$\left( iii\right) $ $\mathcal{M}$ \emph{is uniformly calm at} $\overline{y}%
; $

$\left( iv\right) $ $\mathcal{M}$ \emph{is calm at every} $\left( \overline{y%
},x\right) \in \gph\mathcal{M}.$
\end{rem}

The next three examples show that all converse implications in the previous
remark may fail for a suitable multifunction.

\begin{exa}
\label{Exa_1}\emph{Let }$\mathcal{M}:\mathbb{R}\rightrightarrows \mathbb{R}$%
\emph{\ be given by }$\mathcal{M}\left( y\right) =\left\{ h_{r}\left(
y\right) ,\;r\in \mathbb{N}\right\} ,$ \emph{where} 
\begin{equation*}
h_{r}\left( y\right) =\left\{ 
\begin{array}{lcl}
r+y & \text{\emph{if}} & y\leq \frac{1}{r},\smallskip \\ 
r+\frac{1}{r}+r\left( y-\frac{1}{r}\right) & \text{\emph{if}} & y>\frac{1}{r}%
.%
\end{array}%
\right.
\end{equation*}%
\emph{For }$\overline{y}=0$\emph{, it is easy to check that} $\clm\mathcal{M}%
\left( \overline{y},x\right) =1$ \emph{for all }$x\in \mathcal{M}\left( 
\overline{y}\right) $\emph{. Hence, }$\sup_{x\in \mathcal{M}\left( \overline{%
y}\right) }\clm\mathcal{M}\left( \overline{y},x\right) =1$. \emph{%
Nevertheless, it is impossible to find} $\varepsilon >0$\emph{\ that meets
the conditions for uniform calmness; i.e., }$\uclm\mathcal{M}\left( 
\overline{y}\right) =+\infty $. \emph{More specifically, take }$\varepsilon
_{r}:=r^{-1}+r^{-1/2}$\emph{\ for all }$r\in \mathbb{N},$\ $r\geq 8$\emph{\
(to ensure }$\varepsilon _{r}<1/2),$ \emph{and consider }$%
y_{r}:=r^{-1}+r^{-3/2}$ \emph{and} $x_{r}:=h_{r}\left( y_{r}\right)
=r+r^{-1}+r^{-1/2}\in \mathcal{M}_{\varepsilon _{r}}\left( y_{r}\right) .$ 
\emph{Then}%
\begin{equation*}
\frac{d\left( x_{r},\mathcal{M}\left( 0\right) \right) }{d\left(
y_{r},0\right) }=\frac{r^{-1}+r^{-1/2}}{r^{-1}+r^{-3/2}}\rightarrow +\infty 
\text{ }\emph{as}\text{ }r\rightarrow +\infty .
\end{equation*}
\end{exa}

\begin{exa}
\label{Exa not unif clm}\emph{Consider} $\mathcal{M}:\mathbb{R}%
\longrightarrow \mathbb{R}$ \emph{(single-valued) given by }$\mathcal{M}%
(y)=0 $ \emph{if }$y\leq 0\;$\emph{and} $\mathcal{M}(y)=1$ \emph{if }$y>0.$ 
\emph{It is clear that }$\mathcal{M}$ \emph{is uniformly calm at} $\overline{%
y}=0$\emph{\ (take }$\varepsilon =1/2$ \emph{) but not Lipschitz upper
semicontinuous by just considering }$y_{r}=1/r$\emph{\ for} $r\in \mathbb{N}$%
.
\end{exa}

\begin{exa}
\emph{Let }$\mathcal{M}:\mathbb{R}\rightrightarrows \mathbb{R}$\emph{\ be
given by}%
\begin{equation*}
\mathcal{M}(y)=[0,1]\text{ \emph{if }}y<0,\;\mathcal{M}(y)=[0,+\infty
\lbrack \text{ \emph{if} }y\geq 0.
\end{equation*}%
\emph{It is clear that }$\mathcal{M}$ \emph{is Lipschitz upper
semicontinuous, with zero modulus, at any }$y\in \mathbb{R}$\emph{.
Nevertheless, it is not Hoffman stable at any }$\overline{y}<0.$
\end{exa}

The next theorem establishes that all inequalities in (\ref{eq_inequalities
chain}) become equalities under the convexity of $\gph\mathcal{M}$ together
with the closedness of $\mathcal{M}\left( \overline{y}\right) ,$ provided
that\textbf{\ }$Y$ is a normed space and\textbf{\ }$X$ is a reflexive Banach
space. As an obvious consequence, all properties in Remark \ref{Rem
implication chain} become equivalent in such a case. Firstly, we include two
lemmas.

\begin{lem}
\label{lem 1} Let $X$ be a normed space and $\emptyset \neq C\subset X$ be a
closed set. Take any $x\in X$ and assume that there exists a best
approximation, $\overline{x}$, of $x$ in $C.$ Then $\overline{x}$ is a best
approximation of $x_{\lambda }:=\left( 1-\lambda \right) \overline{x}%
+\lambda x$ in $C$ for all $\lambda \in \lbrack 0,1]$.
\end{lem}

\begin{dem}
Reasoning by contradiction\textbf{, }suppose that for some $\lambda \in
\lbrack 0,1]$ there exists $\hat{x}\in C$ such that $\left\Vert \hat{x}%
-x_{\lambda }\right\Vert <\left\Vert \overline{x}-x_{\lambda }\right\Vert $.
Then 
\begin{align*}
\left\Vert \hat{x}-x\right\Vert & \leq \left\Vert \hat{x}-x_{\lambda
}\right\Vert +\left\Vert x_{\lambda }-x\right\Vert <\left\Vert \overline{x}%
-x_{\lambda }\right\Vert +\left\Vert x_{\lambda }-x\right\Vert \\
& =\lambda \left\Vert \overline{x}-x\right\Vert +\left( 1-\lambda \right)
\left\Vert \overline{x}-x\right\Vert =\left\Vert \overline{x}-x\right\Vert ,
\end{align*}%
which contradicts the fact that $\overline{x}$ is a best approximation of $x$
in $C$.
\end{dem}

In the next result $X$ is assumed to be a reflexive Banach space in order to
ensure the existence of best approximations on nonempty closed convex sets;
see e.g. \cite[Theorem 3.8.1]{ZA02}.

\begin{lem}
\label{lem 2} Let $\mathcal{M}:Y\rightrightarrows X$ be a multifunction
between a normed space $Y$ and a reflexive Banach space $X,$ and assume that 
$\gph\mathcal{M}$ is a nonempty convex set. Let $\overline{y}\in \dom%
\mathcal{M}$ and suppose that $\mathcal{M}\left( \overline{y}\right) $ is
closed. Consider any $\left( y,x\right) \in \gph\mathcal{M}$ and let $%
\overline{x}$ be a best approximation of $x$ in $\mathcal{M}\left( \overline{%
y}\right) $, then 
\begin{equation*}
\frac{d\left( x,\mathcal{M}\left( \overline{y}\right) \right) }{d\left( y,%
\overline{y}\right) }\leq \clm\mathcal{M}\left( \overline{y},\overline{x}%
\right) .
\end{equation*}
\end{lem}

\begin{dem}
By the convexity assumption, for each $\lambda \in \left[ 0,1\right] ,$%
\begin{equation*}
\left( y_{\lambda },x_{\lambda }\right) :=\left( 1-\lambda \right) \left( 
\overline{y},\overline{x}\right) +\lambda \left( y,x\right) \in \gph\mathcal{%
M}.
\end{equation*}%
According to lemma \ref{lem 1}, $\overline{x}$ is also a best approximation
of $x_{\lambda }$ in $\mathcal{M}\left( \overline{y}\right) $, for each $%
\lambda \in \lbrack 0,1]$. Therefore,\textbf{\ } 
\begin{equation*}
\frac{d\left( x,\mathcal{M}\left( \overline{y}\right) \right) }{d\left( y,%
\overline{y}\right) }=\frac{\left\Vert x-\overline{x}\right\Vert }{%
\left\Vert y-\overline{y}\right\Vert }=\frac{\left\Vert x_{\lambda }-%
\overline{x}\right\Vert }{\left\Vert y_{\lambda }-\overline{y}\right\Vert }=%
\frac{d\left( x_{\lambda },\mathcal{M}\left( \overline{y}\right) \right) }{%
d\left( y_{\lambda },\overline{y}\right) },\text{ for all }\lambda \in 
\mathbf{]}0,1\mathbf{]}.
\end{equation*}

Since, letting $\lambda \rightarrow 0,$ we have $\left( y_{\lambda
},x_{\lambda }\right) \rightarrow \left( \overline{y},\overline{x}\right) ,$
by the definition of the calmness modulus (recall (\ref{eq_hof_calm_def}))
we conclude 
\begin{equation*}
\clm\mathcal{M}\left( \overline{y},\overline{x}\right) \geq \limsup_{\lambda
\rightarrow 0}\frac{d\left( x_{\lambda },\mathcal{M}\left( \overline{y}%
\right) \right) }{d\left( y_{\lambda },\overline{y}\right) }=\frac{d\left( x,%
\mathcal{M}\left( \overline{y}\right) \right) }{d\left( y,\overline{y}%
\right) }.
\end{equation*}
\end{dem}

\begin{theo}
\label{The_calmness_Hoffman_convex}Let $\mathcal{M}:Y\rightrightarrows X,$
with $Y$ being a normed space and $X$ being a reflexive Banach space, and
assume that $\gph\mathcal{M}$ is a nonempty convex set. Let $\overline{y}\in %
\dom\mathcal{M}$ with $\mathcal{M}\left( \overline{y}\right) $ closed. Then
one has%
\begin{equation*}
\sup_{x\in \mathcal{M}(\overline{y})}\clm\mathcal{M}\left( \overline{y}%
,x\right) =\uclm\mathcal{M}\left( \overline{y}\right) =\Lipusc\mathcal{M}%
\left( \overline{y}\right) =\Hof\mathcal{M}\left( \overline{y}\right) .
\end{equation*}
\end{theo}

\begin{dem}
We only have to prove $\Hof\mathcal{M}\left( \overline{y}\right) \leq
\sup_{x\in \mathcal{M}(\overline{y})}\clm\mathcal{M}\left( \overline{y}%
,x\right) ,$ according to (\ref{eq_inequalities chain}).

Take any $\left( \widetilde{y},\widetilde{x}\right) \in \gph\mathcal{M}$ and
let $\overline{x}$ be a best approximation of $\widetilde{x}$ in $\mathcal{M}%
(\overline{y})$. Lemma \ref{lem 2} ensures that 
\begin{equation*}
\frac{d\left( \widetilde{x},\mathcal{M}(\overline{y})\right) }{d\left( 
\widetilde{y},\overline{y}\right) }\leq \clm\mathcal{M}\left( \overline{y},%
\overline{x}\right) \leq \sup_{x\in \mathcal{M}(\overline{y})}\clm\mathcal{M}%
\left( \overline{y},x\right) .
\end{equation*}%
Then, recalling (\ref{eq_hof_calm_def}), we conclude 
\begin{equation*}
\Hof\mathcal{M}\left( \overline{y}\right) =\sup_{\left( \widetilde{y},%
\widetilde{x}\right) \in \gph\mathcal{M}}\frac{d(\widetilde{x},\mathcal{M}(%
\overline{y}))}{d\left( \widetilde{y},\overline{y}\right) }\leq \sup_{x\in 
\mathcal{M}(\overline{y})}\clm\mathcal{M}\left( \overline{y},x\right) .
\end{equation*}
\end{dem}

We finish this section by observing that the global Hoffman constant for the
whole graph can be larger than the Hoffman modulus for a specific $\overline{%
y}$. Just consider $\mathcal{M}:\mathbb{R}\rightrightarrows \mathbb{R}$\emph{%
\ }given by 
\begin{equation*}
\mathcal{M}(y)=\left] -\infty ,y\right] \text{ if }\emph{\ }y<0,\;\mathcal{M}%
(y)=\left] -\infty ,0\right] \text{ if }\emph{\ }y\geq 0.
\end{equation*}%
Then clearly $\Hof\mathcal{M}\left( \overline{y}\right) =1$ if$\emph{\ }%
\overline{y}<0\;$and $\Hof\mathcal{M}\left( \overline{y}\right) =0$ if $%
\emph{\ }\overline{y}\geq 0;$ so that $\Hof\mathcal{M}=1.$

\section{Hoffman and calmness moduli for linear semi-infinite inequality
systems}

This section aims to obtain expressions for $\Hof\mathcal{F}$ and $\Hof%
\mathcal{F}\left( \overline{b}\right) ,$ $\overline{b}\in \dom\mathcal{F}$,
in terms of the system's data$.$ These expressions are established in
Theorems \ref{Th_hof:global} and \ref{thm 1}, respectively. The first result
generalizes Theorem \ref{TH_PVZ_const} to the current semi-infinite
framework, while the second provides an alternative expression to Theorem %
\ref{Th_Aze_Cor_subdif} $\left( i\right) ,$ via points inside $\mathcal{F}%
\left( \overline{b}\right) ,$ for locally polyhedral systems$.$ In the case
of finite linear systems Theorem \ref{thm 1} is particularly useful as far
as it establishes an implementable procedure for computing $\Hof\mathcal{F}%
\left( \overline{b}\right) .$

\begin{theo}
\label{Th_hof:global}Consider $\mathcal{F}:C\left( T,\mathbb{R}\right)
\rightrightarrows \mathbb{R}^{n}$ defined in \emph{(\ref{eq_feasible})}. We
have%
\begin{equation*}
\Hof\mathcal{F=}\sup_{_{\substack{ J\subset T\text{ compact}  \\ 0_{n}\notin %
\conv\left\{ a_{t},~t\in J\right\} }}}d_{\ast }\left( 0_{n},\conv\left\{
a_{t},~t\in J\right\} \right) ^{-1}.
\end{equation*}
\end{theo}

\begin{dem}
It is clear that $\Hof\mathcal{F}=\sup_{b\in \dom\mathcal{F}}\Hof\mathcal{F}%
\left( b\right) ,$ and applying Theorem \ref{Th_Aze_Cor_subdif} we have%
\begin{equation}
\Hof\mathcal{F}=\sup_{b\in \dom\mathcal{F}}\sup_{x\notin \mathcal{F}\left(
b\right) }d_{\ast }\left( 0_{n},\mathrm{conv}\left\{ a_{t},\text{ }t\in
J_{b}\left( x\right) \right\} \right) ^{-1}.  \label{eq_hofF}
\end{equation}%
Hence, inequality `$\leq $' comes from (\ref{eq_hofF}) taking into account
that $b\in \dom\mathcal{F}$ and $x\notin \mathcal{F}\left( b\right) $ imply $%
0_{n}\notin \conv\left\{ a_{t},~t\in J_{b}\left( x\right) \right\} $ (recall
Remark \ref{Rem1})$.$ Take also into account that each $J_{b}\left( x\right) 
$ is compact since it is closed in $T$ as far as $J_{b}\left( x\right) $ is
the preimage of $\left\{ f_{b}\left( x\right) \right\} $ by the continuous
function $t\mapsto a_{t}^{\prime }x-b_{t}.$

Let us prove the converse inequality `$\geq $'. Observe that for $%
J=\emptyset $ we have $d_{\ast }\left( 0_{n},\conv\left\{ a_{t},~t\in
J\right\} \right) ^{-1}=d_{\ast }\left( 0_{n},\emptyset \right) ^{-1}=0.$
Fix a nonempty compact set $\widehat{J}\subset T$ such that $0_{n}\notin %
\conv\left\{ a_{t},~t\in \widehat{J}\right\} $ and let us define $\widehat{b}%
\in C\left( T,\mathbb{R}\right) $ such that 
\begin{equation*}
\widehat{J}=J_{\widehat{b}}\left( \widehat{x}\right) ,\text{ for some }%
\widehat{x}\notin \mathcal{F}\left( \widehat{b}\right) ,\text{ }\widehat{b}%
\in \dom\mathcal{F}.
\end{equation*}%
First, by separation, since $0_{n}\notin \conv\left\{ a_{t},~t\in \widehat{J}%
\right\} ,$ there exists $0_{n}\neq \widehat{x}\in \mathbb{R}^{n},$ such
that 
\begin{equation*}
a_{t}^{\prime }\widehat{x}\geq \widehat{x}^{\prime }\widehat{x},\text{ for
all }~t\in \widehat{J},
\end{equation*}%
where $\widehat{x}$ is the best approximation of $0_{n}$ in the compact set $%
\conv\left\{ a_{t},~t\in \widehat{J}\right\} $ with respect to the Euclidean
norm in $\mathbb{R}^{n}.$ Define%
\begin{equation*}
\widehat{b}_{t}:=\max \{a_{t}^{\prime }\widehat{x},\tfrac{1}{2}\widehat{x}%
^{\prime }\widehat{x}\}-\varphi \left( t\right) \tfrac{1}{2}\widehat{x}%
^{\prime }\widehat{x},\text{ }t\in T,
\end{equation*}%
where 
\begin{equation*}
\varphi \left( t\right) =1-d\left( t,\widehat{J}\right) ,\text{ for all }%
t\in T.
\end{equation*}%
Observe that $\widehat{b}\in \dom\mathcal{F}$ since $\widehat{b}_{t}\geq 
\frac{1}{2}\left( 1-\varphi \left( t\right) \right) \widehat{x}^{\prime }%
\widehat{x}\geq 0$ for all $t\in T$ and for instance $0_{n}\in \mathcal{F}%
\left( \widehat{b}\right) .$ On the other hand, $\widehat{x}\notin \mathcal{F%
}\left( \widehat{b}\right) $ since, 
\begin{equation*}
a_{t}^{\prime }\widehat{x}-\widehat{b}_{t}=a_{t}^{\prime }\widehat{x}-\left(
a_{t}^{\prime }\widehat{x}-\varphi \left( t\right) \tfrac{1}{2}\widehat{x}%
^{\prime }\widehat{x}\right) =\tfrac{1}{2}\widehat{x}^{\prime }\widehat{x}>0,%
\text{ if }t\in \widehat{J}.
\end{equation*}%
Finally, observe that%
\begin{equation*}
a_{t}^{\prime }\widehat{x}-\widehat{b}_{t}\leq a_{t}^{\prime }\widehat{x}%
-a_{t}^{\prime }\widehat{x}+\varphi \left( t\right) \tfrac{1}{2}\widehat{x}%
^{\prime }\widehat{x}<\tfrac{1}{2}\widehat{x}^{\prime }\widehat{x},\text{
whenever }t\in T\setminus \widehat{J}.
\end{equation*}%
So,%
\begin{equation*}
\widehat{J}=\left\{ t\in T\mid a_{t}^{\prime }\widehat{x}-\widehat{b}_{t}=f_{%
\widehat{b}}\left( \widehat{x}\right) \right\} ,
\end{equation*}%
in other words, $\widehat{J}=J_{\widehat{b}}\left( \widehat{x}\right) ,$
which finishes the proof.
\end{dem}

\begin{rem}
\emph{Theorem \ref{Th_hof:global} is the only result in this paper which
uses the fact that }$T$\emph{\ is assumed to be a compact} metric\emph{\
space. The rest of results work for }$T$\emph{\ being a compact Hausdorff
space, which is the framework of the so-called }continuous systems\emph{\ in 
\cite{libro}.}
\end{rem}

The rest of this section is focussed on $\Hof\mathcal{F}\left( \overline{b}%
\right) \,,$ provided that $\overline{b}\in \dom\mathcal{F}$. To start with,
as a consequence of Theorem \ref{The_calmness_Hoffman_convex}, we always
have 
\begin{equation}
\Hof\mathcal{F}\left( \overline{b}\right) =\sup_{x\in \mathcal{F}\left( 
\overline{b}\right) }\clm\mathcal{F}\left( \overline{b},x\right) =\sup_{x\in 
\mathrm{bd}\mathcal{F}\left( \overline{b}\right) }\clm\mathcal{F}\left( 
\overline{b},x\right) ,\text{ }\overline{b}\in \dom\mathcal{F},
\label{eq_h_Cal}
\end{equation}%
where the last equality comes from the fact that $\clm\mathcal{F}\left( 
\overline{b},x\right) =0$ when $x\in \mathrm{int}\mathcal{F}\left( \overline{%
b}\right) $ (the trivial case $\mathrm{bd}\mathcal{F}\left( \overline{b}%
\right) =\emptyset $, equivalently $\mathcal{F}\left( \overline{b}\right) =%
\mathbb{R}^{n},$ is included; recall $\sup \emptyset :=0).$ From now on we
are devoted to refine (\ref{eq_h_Cal}) by replacing $\mathrm{bd}\mathcal{F}%
\left( \overline{b}\right) $ with a smaller subset. The concluding result is
Theorem \ref{thm 1}. First, we establish some technical results.

\begin{prop}
\label{Prop_2} Let $x^{1},x^{2}\in \mathrm{bd}\mathcal{F}\left( \overline{b}%
\right) $ such that $T\left( x^{1}\right) \subset T\left( x^{2}\right) $.
Then,

$\left( i\right) $ $\mathrm{end}\partial f_{\overline{b}}\left( x^{1}\right)
\subset \mathrm{end}\partial f_{\overline{b}}\left( x^{2}\right) $;

$\left( ii\right) $ If the regularity condition (\ref%
{eq_regularity_condition}) is held at $x^{i},$ $i=1,2,$ then 
\begin{equation*}
\clm\mathcal{F}\left( \overline{b},x^{1}\right) \leq \clm\mathcal{F}\left( 
\overline{b},x^{2}\right) .
\end{equation*}
\end{prop}

\begin{dem}
$\left( i\right) $ First, $x^{i}\in \mathrm{bd}\mathcal{F}\left( \overline{b}%
\right) $ implies $f_{\overline{b}}\left( x^{i}\right) =0,$ and so $T\left(
x^{i}\right) \neq \emptyset ,$ $i=1,2,$ by the compactness of $T$ together
with the continuity of $t\mapsto \tbinom{a_{t}}{\overline{b}_{t}}.$ Recall
that, $\partial f_{\overline{b}}\left( x^{i}\right) =\mathrm{conv}\left\{
a_{i},~i\in T\left( x^{i}\right) \right\} ,$ $i=1,2,$ hence $\partial f_{%
\overline{b}}\left( x^{1}\right) \subset \partial f_{\overline{b}}\left(
x^{2}\right) .$

Assume, arguing by contradiction, that there exists $a\in \mathrm{end}%
\partial f_{\overline{b}}\left( x^{1}\right) \setminus \mathrm{end}\partial
f_{\overline{b}}\left( x^{2}\right) .$ Since, by compactness, $\mathrm{end}%
\partial f_{\overline{b}}\left( x^{1}\right) \subset \partial f_{\overline{b}%
}\left( x^{1}\right) \subset \partial f_{\overline{b}}\left( x^{2}\right) ,$
we have $a\in \partial f_{\overline{b}}\left( x^{2}\right) \setminus \mathrm{%
end}\partial f_{\overline{b}}\left( x^{2}\right) .$ Then we have $\lambda
a\in \partial f_{\overline{b}}\left( x^{2}\right) $ for some $\lambda >1$
and we can write%
\begin{equation}
\lambda a=\sum\limits_{t\in T\left( x^{1}\right) }\lambda
_{t}a_{t}+\sum\limits_{t\in T\left( x^{2}\right) \backslash T\left(
x^{1}\right) }\lambda _{t}a_{t},  \label{eq_00000}
\end{equation}%
for some $\left\{ \lambda _{t}\right\} _{t\in T\left( x^{2}\right) }\subset 
\mathbb{R}_{+}$ such that $\left\{ \lambda _{t}\mid \lambda _{t}\neq 0,\text{
}t\in T\left( x^{2}\right) \right\} $ is a finite set.

On the other hand, consider\textbf{\ }$d:=x^{1}-x^{2}$ and observe that, 
\begin{equation*}
\left\{ 
\begin{array}{l}
a_{t}^{\prime }d=0,\text{ }t\in T\left( x^{1}\right) , \\ 
a_{t}^{\prime }d=a_{t}^{\prime }x^{1}-a_{t}^{\prime }x^{2}<\overline{b}_{t}-%
\overline{b}_{t}=0,\text{ }t\in T\left( x^{2}\right) \backslash T\left(
x^{1}\right) .%
\end{array}%
\right.
\end{equation*}%
Then, multiplying (with the inner product) both members of (\ref{eq_00000})
by $d$, we deduce%
\begin{equation*}
0=\lambda a^{\prime }d=\sum\limits_{t\in T\left( x^{2}\right) \backslash
T\left( x^{1}\right) }\lambda _{t}a_{t}^{\prime }d,
\end{equation*}%
which yields $\lambda _{t}=0$ for all $t\in T\left( x^{2}\right) \backslash
T\left( x^{1}\right) .$ So, we attain the contradiction $\lambda
a=\sum\nolimits_{t\in T\left( x^{1}\right) }\lambda _{t}a_{t}\in \partial
f_{\overline{b}}\left( x^{1}\right) .$

Statement $\left( ii\right) $ follows straightforwardly from Theorem \ref%
{Th_LiMengYan}.
\end{dem}

The following example shows that the regularity condition assumed in
statement $\left( ii\right) $ of the previous proposition is not
superfluous. The example comes from modifying Example 1 in \cite{CLPT14}
(revisited in \cite[Example 3.3]{LMY18}).

\begin{exa}
\label{Exa1}\emph{Let us consider the system, in} $\mathbb{R}^{2}$ \emph{%
endowed with the Euclidean norm}$,$ \emph{given by} 
\begin{equation*}
\sigma \left( \overline{b}\right) :=\left\{ 
\begin{tabular}{rl}
$t\left( \cos t\right) x_{1}+t\left( \sin t\right) x_{2}\leq t,$ & $t\in %
\left[ 0,\pi \right] ,$ \\ 
$x_{1}\leq 1,$ & $t=4,$ \\ 
$-x_{1}-x_{2}\leq 1,$ & $t=5$%
\end{tabular}%
\right\} ;
\end{equation*}%
\emph{i.e.,} $T:=\left[ 0,\pi \right] \cup \{4,5\},$ $a_{t}:=t\left( \cos
t,\sin t\right) ^{\prime },$ \emph{for} $t\in \left[ 0,\pi \right] ,$ $%
a_{4}:=\left( 1,0\right) ^{\prime }$ \emph{and }$a_{5}:=\left( -1,-1\right)
^{\prime };$ $\overline{b}\in C\left( \left[ 0,\pi \right] \cup \{4,5\},%
\mathbb{R}\right) $ \emph{is given by }$\overline{b}_{t}=t,$ $t\in \left[
0,\pi \right] $ , $\overline{b}_{4}=1,$\emph{\ and } $\overline{b}_{5}=1.$ 
\emph{\ Consider the feasible points} $x^{1}=\left( 1,0\right) ^{\prime }$ 
\emph{and }$x^{2}=\left( 1,-2\right) ^{\prime }.$

\begin{figure}[h]
	\centering
	\includegraphics[scale=0.4]{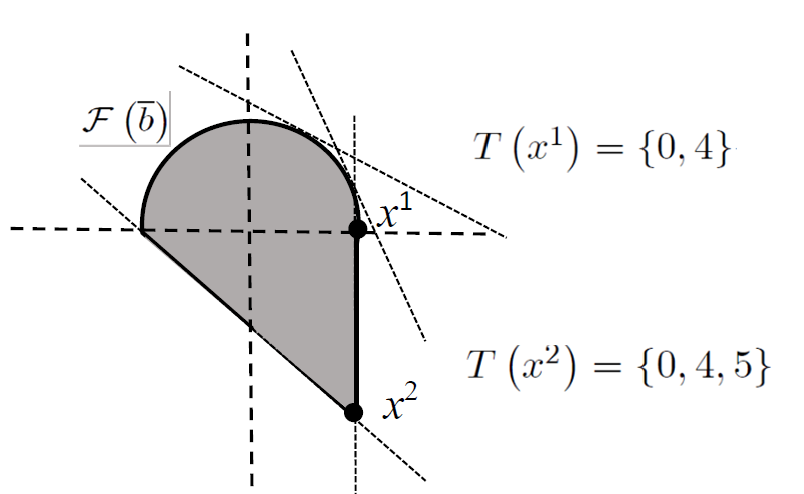} 
	\caption{Illustration of Example \ref{Exa1}}
\end{figure}

\emph{As proved in \cite[Example 1]{CLPT14}, we have that} 
\begin{equation*}
\clm\mathcal{F}\left( \overline{b},x^{1}\right) =+\infty .
\end{equation*}%
\emph{Alternatively, we can apply Theorem \ref{Th_Aze_Cor_subdif}}$\left(
ii\right) $ \emph{with sequence }$x^{r}=\left( 1+\frac{1}{r}\right) \binom{%
\cos \frac{1}{r}}{\sin \frac{1}{r}}.$ \emph{It is clear that the regularity
condition (\ref{eq_regularity_condition}) is not satisfied at }$x^{1}.\,\ $%
\emph{Indeed} $\left( 0,1\right) ^{\prime }\in A\left( x^{1}\right) ^{\circ
}=\mathrm{cone}\left\{ \left( 1,0\right) ^{\prime }\right\} ^{\circ }=%
\mathbb{R}_{-}\times \mathbb{R}$\emph{, but} $x^{1}+\varepsilon \left(
0,1\right) ^{\prime }\notin \mathcal{F}\left( \overline{b}\right) $ \emph{%
for any} $\varepsilon >0.$ \emph{Moreover,} 
\begin{equation*}
\partial f_{\overline{b}}\left( x^{1}\right) =\mathrm{conv}\left\{ \left(
0,0\right) ^{\prime },\left( 1,0\right) ^{\prime }\right\} 
\end{equation*}%
\emph{and} $\mathrm{end}\partial f_{\overline{b}}\left( x^{1}\right)
=\left\{ \left( 1,0\right) ^{\prime }\right\} .$ \emph{Hence, }$\clm\mathcal{%
F}\left( \overline{b},x^{1}\right) \neq d_{\ast }\left( 0_{2},\mathrm{end}%
\partial f_{\overline{b}}\left( x^{1}\right) \right) ^{-1}.$

\emph{With respect to point }$x^{2},$ \emph{one easily sees that condition (%
\ref{eq_regularity_condition}) is satisfied}, \emph{where} $A\left(
x^{2}\right) ^{\circ }=\left\{ u\in \mathbb{R}^{2}\mid -u_{1}-u_{2}\leq
0,u_{1}\leq 0\right\} .\,\ $\emph{In this case,} $\partial f_{\overline{b}%
}\left( x^{2}\right) =\mathrm{conv}\left\{ \left( 0,0\right) ^{\prime
},\left( 1,0\right) ^{\prime },\left( -1,-1\right) ^{\prime }\right\} .\,\ $%
\emph{Hence, from Theorem \ref{Th_LiMengYan} we have} 
\begin{eqnarray*}
\clm\mathcal{F}\left( \overline{b},x^{2}\right)  &=&d_{\ast }\left( 0_{2},%
\mathrm{end}\partial f_{\overline{b}}\left( x^{2}\right) \right) ^{-1} \\
&=&d_{\ast }\left( 0_{2},\mathrm{conv}\left\{ \left( 1,0\right) ^{\prime
},\left( -1,-1\right) ^{\prime }\right\} \right) ^{-1}=\sqrt{5}.
\end{eqnarray*}%
$\bigskip $
\end{exa}

\begin{prop}
\label{Prop3}Let $C$ be a nonempty closed convex subset of $\mathbb{R}^{n}$
different from a singleton with $\mathrm{extr\,}C\neq \emptyset $ and let $%
x^{0}\in C\backslash \mathrm{extr\,}C.$ Then, there exist $y^{0}\in \mathrm{%
extr\,}C,$ $z^{0}\in C,$ and $\mu \in \left] 0,1\right[ $ such that $%
x^{0}=\left( 1-\mu \right) y^{0}+\mu z^{0}.$
\end{prop}

\begin{dem}
The assumption $\mathrm{extr\,}C\neq \emptyset $ is equivalent to the fact
that $C$ contains no lines (i.e., its lineality space is $\left\{
0_{n}\right\} $). According to \cite[Corollary 14.6.1]{Rock70}, this is also
equivalent to $\mathrm{int}\left( O^{+}C\right) ^{\circ }\neq \emptyset ,$
recalling that $O^{+}C$ is the recession cone of $C.$ Pick $0_{n}\neq u\in 
\mathrm{int}\left( O^{+}C\right) ^{\circ }$ and consider 
\begin{equation*}
K:=C\cap \left\{ x\in \mathbb{R}^{n}\mid u^{\prime }x\geq u^{\prime
}x^{0}-1\right\} .
\end{equation*}%
Let us see that $K$ is bounded, i.e., $O^{+}K=\left\{ 0_{n}\right\} $ (see 
\cite[Theorem 8.4]{Rock70}). Reasoning by contradiction, assume the
existence of $0_{n}\neq v\in O^{+}K.$ Then $x^{0}+\lambda v\in K$ and,
accordingly, $u^{\prime }\left( x^{0}+\lambda v\right) \geq u^{\prime
}x^{0}-1$ for all $\lambda >0.$ Letting $\lambda \rightarrow +\infty $ we
obtain $u^{\prime }v\geq 0$. On the other hand, $v\in O^{+}C$ and, for $%
\alpha >0$ small enough, we have $u+\alpha v\in \left( O^{+}C\right) ^{\circ
},$ yielding the contradiction $0\geq \left( u+\alpha v\right) ^{\prime
}v\geq \alpha v^{\prime }v.$

Once we know that $K$ is a nonempty convex compact set, by applying the
Minkowski-Carath\'{e}odory theorem (see, e.g., \cite[Theorem 8.11]{Simon11}%
), we have $K=\mathrm{conv}\left( \mathrm{extr\,}K\right) ,$ and we can
write 
\begin{equation}
x^{0}=\sum_{i=1}^{k}\lambda _{i}x^{i}  \label{eqMC}
\end{equation}%
with $\left\{ x^{1},...,x^{k}\right\} \subset \mathrm{extr\,}K$ being
affinely independent, $\sum_{i=1}^{k}\lambda _{i}=1,$ and $\lambda _{i}>0$
for all \thinspace $i=1,...,k.$ Clearly it is not restrictive to assume $%
u^{\prime }x^{1}\geq u^{\prime }x^{0}$, which easily entails $x^{1}\in 
\mathrm{extr\,}C$. More in detail, if $x^{1}$ were a midpoint of distinct
points in $C,$ we could replace these points with others in the same segment
verifying $u^{\prime }x\geq u^{\prime }x^{0}-1,$ and hence these points
would be in $K,$ contradicting $x^{1}\in \mathrm{extr\,}K.$

On the other hand, by applying \cite[Theorem A.7]{libro}, (\ref{eqMC})
entails that $x^{0}$ is in the relative interior of $\mathrm{conv}\left\{
x^{1},...,x^{k}\right\} $ (i.e., the interior relative to the affine hull of
these points), and then $z^{0}:=x^{1}+\beta \left( x^{0}-x^{1}\right) \in 
\mathrm{conv}\left\{ x^{1},...,x^{k}\right\} \subset C$ for a small enough $%
\beta >1.$ Finally, let us write 
\begin{equation*}
x^{0}=\left( 1-\tfrac{1}{\beta }\right) x^{1}+\tfrac{1}{\beta }z^{0},
\end{equation*}%
which provides the aimed result with $y^{0}=x^{1}$ and $\mu =\tfrac{1}{\beta 
}.$
\end{dem}

The following theorem appeals to locally polyhedral (LOP, in brief) systems.
Recall that given $\overline{b}\in \dom\mathcal{F},$ $\sigma \left( 
\overline{b}\right) $ is a LOP system iff 
\begin{equation}
D\left( \mathcal{F}\left( \overline{b}\right) ,\overline{x}\right) =A\left( 
\overline{x}\right) ^{\circ },\text{ for all }\overline{x}\in \mathcal{F}%
\left( \overline{b}\right) ,  \label{eq_LOP}
\end{equation}%
where $D\left( \mathcal{F}\left( \overline{b}\right) ,\overline{x}\right) $
denotes the \emph{cone of feasible directions} of $\mathcal{F}\left( 
\overline{b}\right) $ at $\overline{x};$ i.e., $d\in D\left( \mathcal{F}%
\left( \overline{b}\right) ,\overline{x}\right) $ if there exists $%
\varepsilon >0$ such that $\overline{x}+\alpha d\in \mathcal{F}\left( 
\overline{b}\right) $ for all $\alpha \in \left[ 0,\varepsilon \right] .$
See \cite{AGL98} for a comprehensive analysis of LOP systems (see also \cite%
{libro}). At this moment we recall a characterization of LOP systems in
terms of the regularity condition (\ref{eq_LOP}) which can be derived from
Corollary 3.3 in \cite{LMY18}.

\begin{lem}
\emph{(see \cite[Corollary 3.3]{LMY18})} \label{Lem1}Let $\overline{b}\in %
\dom\mathcal{F}$. The following conditions are equivalent:

$\left( i\right) $ $D\left( \mathcal{F}\left( \overline{b}\right) ,\overline{%
x}\right) =A\left( \overline{x}\right) ^{\circ },$ for all $\overline{x}\in 
\mathcal{F}\left( \overline{b}\right) ,$

$\left( ii\right) $ The regularity condition (\ref{eq_regularity_condition})
is held at any $\overline{x}\in \mathcal{F}\left( \overline{b}\right) .$
\end{lem}

From now on we consider the set 
\begin{equation}
\mathcal{E}\left( \overline{b}\right) :=\extr\left( \mathcal{F}\left( 
\overline{b}\right) \cap \spann\left\{ a_{t},~t\in T\right\} \right) ,\text{
with }\overline{b}\in \dom\mathcal{F}.  \label{eq_Eb}
\end{equation}%
Observe that, $\mathcal{E}\left( \overline{b}\right) $ is always a nonempty
and finite set when $T$ is finite; moreover,%
\begin{equation*}
\mathcal{E}\left( \overline{b}\right) =\extr\mathcal{F}\left( \overline{b}%
\right) \Leftrightarrow \extr\mathcal{F}\left( \overline{b}\right) \neq
\emptyset ;
\end{equation*}%
in fact, $\extr\mathcal{F}\left( \overline{b}\right) \neq \emptyset $ if and
only if $\mathcal{F}\left( \overline{b}\right) $ does not contain any line,
which is equivalent to the fact that $\spann\left\{ a_{t},~t\in T\right\} =%
\mathbb{R}^{n}.\ $This construction is inspired by the one of \cite[p. 142]%
{Li94}, and used in \cite{GCLT18} to compute the calmness modulus of the
optimal value function of finite linear optimization problems.

\begin{theo}
\label{thm 1} Let $\overline{b}\in \dom\mathcal{F}$ and assume that $\sigma
\left( \overline{b}\right) $ is a LOP system$.$ Then 
\begin{equation*}
\Hof\mathcal{F}\left( \overline{b}\right) =\sup_{x\in \mathcal{E}\left( 
\overline{b}\right) }\clm\mathcal{F}\left( \overline{b},x\right) =\sup_{x\in 
\mathcal{E}\left( \overline{b}\right) }\sup_{D\in \mathcal{D}\left( x\right)
}d_{\ast }\left( 0_{n},\mathrm{conv}\left\{ a_{t},\text{ }t\in D\right\}
\right) ^{-1}.
\end{equation*}
\end{theo}

\begin{dem}
To start with, we recall equation (\ref{eq_h_Cal}): 
\begin{equation*}
\Hof\mathcal{F}\left( \overline{b}\right) =\sup_{x\in \mathrm{bd}\mathcal{F}%
\left( \overline{b}\right) }\clm\mathcal{F}\left( \overline{b},x\right) .
\end{equation*}%
Since $\mathcal{E}\left( \overline{b}\right) \subset \mathrm{bd}\mathcal{F}%
\left( \overline{b}\right) ,$ the inequality $\Hof\mathcal{F}\left( 
\overline{b}\right) \geq \sup_{x\in \mathcal{E}\left( \overline{b}\right) }%
\clm\mathcal{F}\left( \overline{b},x\right) $ follows trivially.

Let us see that $\Hof\mathcal{F}\left( \overline{b}\right) \leq \sup_{x\in 
\mathcal{E}\left( \overline{b}\right) }\clm\mathcal{F}\left( \overline{b}%
,x\right) $. Specifically, let us prove that for every $x\in \mathrm{bd}%
\mathcal{F}\left( \overline{b}\right) $ there exists $\widetilde{x}\in 
\mathcal{E}\left( \overline{b}\right) $ such that $\clm\mathcal{F}\left( 
\overline{b},x\right) \leq \clm\mathcal{F}\left( \overline{b},\widetilde{x}%
\right) .$

Fix arbitrarily $x\in \mathrm{bd}\mathcal{F}\left( \overline{b}\right) $ and
write $x=y+z$, where $y\in \spann\left\{ a_{t},~t\in T\right\} $ and $z\in
\left\{ a_{t},~t\in T\right\} ^{\perp }$ (the orthogonal subspace to $%
\left\{ a_{t},~t\in T\right\} $). Since $a_{t}^{\prime }x=a_{t}^{\prime }y$
for all $t\in T,$ $y\in \mathrm{bd}\mathcal{F}\left( \overline{b}\right) $
and 
\begin{equation*}
T(x)=T(y).
\end{equation*}%
Hence, applying Proposition \ref{Prop_2}$\left( ii\right) $ (recall Lemma %
\ref{Lem1})$,$ we have 
\begin{equation}
\clm\mathcal{F}\left( \overline{b},x\right) =\clm\mathcal{F}\left( \overline{%
b},y\right) .  \label{eq_111}
\end{equation}

Let us denote 
\begin{equation*}
C=\mathcal{F}\left( \overline{b}\right) \cap \spann\left\{ a_{t},~t\in
T\right\} ,
\end{equation*}%
which satisfies $\mathrm{extr\,}C\neq \emptyset .$ If $y\in \mathrm{extr\,}C=%
\mathcal{E}\left( \overline{b}\right) ,$ we are done. Otherwise, if $y\in
C\backslash \mathrm{extr\,}C,$ we can apply Proposition \ref{Prop3} and
conclude the existence of $\widetilde{x}\in \mathrm{extr\,}C,$ $\widetilde{z}%
\in C,$ and $\mu \in \left] 0,1\right[ $ such that $y=\left( 1-\mu \right) 
\widetilde{x}+\mu \widetilde{z}.$ Observe that%
\begin{equation*}
T(y)\subset T(\widetilde{x}),
\end{equation*}%
since $a_{t}^{\prime }y=b_{t}$ implies $\left( 1-\mu \right) a_{t}^{\prime }%
\widetilde{x}+\mu a_{t}^{\prime }\widetilde{z}=b_{t},$ which entails $%
a_{t}^{\prime }\widetilde{x}=a_{t}^{\prime }\widetilde{z}=b_{t}$ (because
both $\widetilde{x},\widetilde{z}\in \mathcal{F}\left( \overline{b}\right)
). $ So, we conclude the aimed inequality 
\begin{equation*}
\clm\mathcal{F}\left( \overline{b},y\right) \leq \clm\mathcal{F}\left( 
\overline{b},\widetilde{x}\right) ,
\end{equation*}%
which together with (\ref{eq_111}) yields%
\begin{equation*}
\clm\mathcal{F}\left( \overline{b},x\right) \leq \clm\mathcal{F}\left( 
\overline{b},\widetilde{x}\right) ,\text{ with }\widetilde{x}\in \mathcal{E}%
\left( \overline{b}\right) .
\end{equation*}
\end{dem}

\subsection{On the finite case}

This subsection gathers some specifics on finite linear systems. Thus, along
this subsection, we assume that $T$\ is finite, in which case, for a fixed $%
\left( \overline{b},\overline{x}\right) \in \gph\mathcal{F},$ $\mathcal{D}%
\left( \overline{x}\right) $ is also finite and, clearly 
\begin{equation*}
\cup _{D\in \mathcal{D}\left( \overline{x}\right) }\mathrm{conv}\left\{
a_{t},\text{ }t\in D\right\} =\mathrm{end}\partial f\left( \overline{x}%
\right) 
\end{equation*}%
is a closed set; moreover, $\mathcal{E}\left( \overline{b}\right) $ is also
finite and $\clm\mathcal{F}\left( \overline{b},\overline{x}\right) $ and $%
\Hof\mathcal{F}\left( \overline{b}\right) $ can be computed through the
implementable computations: 
\begin{eqnarray*}
\clm\mathcal{F}\left( \overline{b},\overline{x}\right)  &=&\max_{D\in 
\mathcal{D}\left( \overline{x}\right) }d_{\ast }\left( 0_{n},\mathrm{conv}%
\left\{ a_{t},\text{ }t\in D\right\} \right) ^{-1}. \\
\Hof\mathcal{F}\left( \overline{b}\right)  &=&\max_{x\in \mathcal{E}\left( 
\overline{b}\right) }\clm\mathcal{F}\left( \overline{b},x\right) .
\end{eqnarray*}

In addition, as a consequence of Theorem \ref{thm 1}, we can write 
\begin{equation*}
\Hof\mathcal{F}=\max_{b\in \dom\mathcal{F}}\Hof\mathcal{F}\left( b\right)
=\max_{b\in \dom\mathcal{F}}\max_{x\in \mathcal{E}\left( b\right) }\clm%
\mathcal{F}\left( b,x\right) .
\end{equation*}%
Indeed, if the maximum in (\ref{eq_2}) in Theorem \ref{TH_PVZ_const} is
attained at $J\subset T$ such that\emph{\ }$\rank A_{J}=\rank A$\emph{\ }and%
\emph{\ }$\{a_{t},~t\in J\}$\emph{\ }is linearly independent, we have\emph{\ 
}%
\begin{equation*}
\Hof\mathcal{F}=\Hof\mathcal{F}\left( b^{J}\right) =\clm\mathcal{F}\left(
b^{J},0_{n}\right) ,
\end{equation*}%
where $b^{J}$\emph{\ }is defined as\emph{\ }$b_{t}^{J}=0$\ if $t\in J$\emph{%
\ }and $b_{t}^{J}=1$\emph{\ }otherwise.

Finally, we observe that Proposition \ref{Prop_2} $\left( i\right) $ admits
a refinement in this finite case, which is written in the following result.

\begin{prop}
\label{prop4} Let $x^{1},x^{2}\in \mathrm{bd}\mathcal{F}\left( \overline{b}%
\right) $ such that $T\left( x^{1}\right) \subset T\left( x^{2}\right) $.
Then, $\mathcal{D}\left( x^{1}\right) \subset \mathcal{D}\left( x^{2}\right)
.$
\end{prop}

\begin{dem}
Given $D\in \mathcal{D}\left( x^{1}\right) $, let us see that $D\in \mathcal{%
D}\left( x^{2}\right) $. First, consider\textbf{\ }$d:=x^{1}-x^{2}$ and
observe that, 
\begin{equation*}
\left\{ 
\begin{array}{l}
a_{t}^{\prime }d=0,\text{ }t\in T\left( x^{1}\right) , \\ 
a_{t}^{\prime }d=a_{t}^{\prime }x^{1}-a_{t}^{\prime }x^{2}<\overline{b}_{t}-%
\overline{b}_{t}=0,\text{ }t\in T\left( x^{2}\right) \backslash T\left(
x^{1}\right) .%
\end{array}%
\right. 
\end{equation*}%
Now, recalling (\ref{SystemD2}), the fact that\textbf{\ }$D\in \mathcal{D}%
\left( x^{1}\right) $ ensures the existence of $\overline{d}\in \mathbb{R}%
^{n}$ such that 
\begin{equation*}
\left\{ 
\begin{array}{l}
a_{t}^{\prime }\overline{d}=1,\text{ }t\in D, \\ 
a_{t}^{\prime }\overline{d}<1,\text{ }t\in T\left( x^{1}\right) \backslash D.%
\end{array}%
\right. 
\end{equation*}%
For every $\alpha >0,$ we consider a new vector $d_{\alpha }:=\overline{d}%
+\alpha d;$ observe that 
\begin{equation*}
\left\{ 
\begin{array}{l}
a_{t}^{\prime }d_{\alpha }=a_{t}^{\prime }\overline{d}+\alpha a_{t}^{\prime
}d=1,\text{ }t\in D, \\ 
a_{t}^{\prime }d_{\alpha }=a_{t}^{\prime }\left( \overline{d}+\alpha
d\right) <1,\text{ }t\in T\left( x^{1}\right) \backslash D.%
\end{array}%
\right. 
\end{equation*}%
Since $a_{t}^{\prime }d<0$\ for\textbf{\ }$t\in T\left( x^{2}\right)
\backslash T\left( x^{1}\right) ,$ we can choose $\alpha $ large enough (any 
$\alpha >\max_{t\in T\left( x^{2}\right) \backslash T\left( x^{1}\right) }%
\frac{a_{t}^{\prime }\overline{d}-1}{-a_{t}^{\prime }d}$ will do it) to make 
$a_{t}^{\prime }\left( \overline{d}+\alpha d\right) <1$ for all $t\in
T\left( x^{2}\right) \backslash T\left( x^{1}\right) $. This proves $D\in 
\mathcal{D}\left( x^{2}\right) .$
\end{dem}

The following example shows that the previous proposition does not hold in
the semi-infinite framework.

\begin{exa}
\label{exa2}\emph{Let us consider the system, in} $\mathbb{R}^{2}$ \emph{%
endowed with the Euclidean norm}$,$ \emph{given by} 
\begin{equation*}
\sigma \left( \overline{b}\right) :=\left\{ 
\begin{tabular}{rl}
$\left( 1+t\cos t\right) x_{1}+\left( t\sin t\right) x_{2}$ & $\leq 0,$ $%
t\in \left[ 0,\frac{\pi }{2}\right] $%
\end{tabular}%
\right\} ;
\end{equation*}%
\emph{and take }$x^{1}=\left( 0,-1\right) ^{\prime }$ \emph{and }$%
x^{2}=\left( 0,0\right) ^{\prime }.$ \emph{Then }$T\left( x^{1}\right)
=\{0\}\subset \left[ 0,\frac{\pi }{2}\right] =T\left( x^{2}\right) .$ \emph{%
We have}%
\begin{equation*}
\{0\}\in \mathcal{D}\left( x^{1}\right) \setminus \mathcal{D}\left(
x^{2}\right) .
\end{equation*}%
\emph{To check that }$\{0\}\notin \mathcal{D}\left( x^{2}\right) $ \emph{%
observe that the system, in the variable }$d=\left( d_{1},d_{2}\right)
^{\prime }\in \mathbb{R}^{2},$%
\begin{equation*}
\left\{ 
\begin{tabular}{rll}
$d_{1}=1,$ & $\left( 1+t\cos t\right) d_{1}+\left( t\sin t\right) d_{2}<1$ & 
$,$ $t\in \left] 0,\frac{\pi }{2}\right] $%
\end{tabular}%
\right\}
\end{equation*}%
\emph{is inconsistent.}
\end{exa}

\section{Conclusions and perspectives}

We have analyzed different properties oriented to quantify the global,
semi-local and local Hoffman\textbf{\ }behavior of set-valued mappings
between metric spaces, where by `semi-local' we mean the study of the whole
image set with respect to parameter perturbations (a similar use of this
term can be found, for instance, in \cite[Definition 2.1]{YYK08}), yielding
to the known Lipschitz upper semicontinuity when the study is concentrated
around a nominal parameter. Local properties, as calmness, are focussed on
the behavior of the multifunction around a fixed element of its graph. The
corresponding moduli are analyzed. Both Hoffman stability (\ref{eq_Hof
stable_b}) and uniform calmness (\ref{eq uniform calmness}) constitute
intermediate steps between calmness and global Hoffman properties. All these
semi-local properties are shown to be equivalent (and with the same
rate/modulus) for convex-graph multifunctions taking closed values in a
reflexive Banach space (Theorem \ref{The_calmness_Hoffman_convex}). This is
the case of the feasible set mapping, $\mathcal{F},$ associated with a
continuous linear semi-infinite inequality system parameterized with respect
to the right-hand side. At this moment, let us comment that paper \cite%
{CCP22} analyzes the upper Lipschitz behavior of the optimal set mapping, $%
\mathcal{F}^{op},$\ in finite linear programming, which does not have a
convex graph. Appealing to a certain concept of directional convexity
introduced in that paper, \cite{CCP22} establishes a counterpart for the
optimal set mapping of formula\textbf{\ }%
\begin{equation*}
\Lipusc\mathcal{F}(\overline{b})=\sup_{x\in \mathcal{F}(\overline{b})}\clm%
\mathcal{F}\left( \overline{b},x\right) .
\end{equation*}%
However, it is shown there that the Hoffman and Lipschitz upper
semicontinuity moduli do not coincide when applied to $\mathcal{F}^{op}$ at
a nominal parameter.

For this feasible set mapping we succeed in giving the following formula for
the global Hoffman constant (Theorem \ref{Th_hof:global}), which extends to
the current semi-infinite framework some previous results for finite systems,%
\begin{equation*}
\Hof\mathcal{F=}\sup_{_{\substack{ J\subset T\text{ compact} \\ 0_{n}\notin %
\conv\left\{ a_{t},~t\in J\right\} }}}d_{\ast }\left( 0_{n},\conv\left\{
a_{t},~t\in J\right\} \right) ^{-1}.
\end{equation*}%
With respect to the semi-local measure, $\Hof\mathcal{F}\left( \overline{b}%
\right) ,$ when confined to locally polyhedral systems (which includes
finite systems), Theorem \ref{thm 1} provides a point-based formula
involving exclusively some feasible points and the nominal data $a_{t}$'s
and $\overline{b}_{t}$'s:%
\begin{equation}
\Hof\mathcal{F}\left( \overline{b}\right) =\sup_{x\in \mathcal{E}\left( 
\overline{b}\right) }\sup_{D\in \mathcal{D}\left( x\right) }d_{\ast }\left(
0_{n},\mathrm{conv}\left\{ a_{t},\text{ }t\in D\right\} \right) ^{-1},
\label{eq_4}
\end{equation}%
where $\mathcal{E}\left( \overline{b}\right) $ is defined in (\ref{eq_Eb}).
When $T$ is finite (and hence $\mathcal{E}\left( \overline{b}\right) $ and
each $\mathcal{D}\left( x\right) $ also are), the previous expression yields
a specially computable procedure. It provides an alternative approach to the
one given in \cite{AzCo02} via points outside the feasible set:%
\begin{equation*}
\Hof\mathcal{F}(\overline{b})=\sup_{x\neq \mathcal{F}\left( \overline{b}%
\right) }d_{\ast }\left( 0_{n},\mathrm{conv}\left\{ a_{t},\text{ }t\in J_{%
\overline{b}}\left( x\right) \right\} \right) ^{-1}.
\end{equation*}%
The problem of finding an expression for $\Hof\mathcal{F}\left( \overline{b}%
\right) $ in the line of (\ref{eq_4}) for not locally polyhedral systems
remains as open problem. A crucial step here is to extended Theorem \ref%
{Th_LiMengYan} about the calmness modulus (traced out from \cite{LMY18}) to
more general semi-infinite system.

\textbf{Acknowledgement: }The authors are indebted to the anonymous referees
and the Associate Editor for their valuable critical comments, which have
definitely improved the original version of the paper.

\end{document}